\documentclass[a4,11pt]{article}
\usepackage{gastex}
\usepackage{xy}
\usepackage{amsgen}
\usepackage{amsmath}
\usepackage{amstext}
\usepackage{amsbsy}
\usepackage{amsopn}
\usepackage{amsfonts}
\usepackage{amssymb}
\usepackage{eepic}
\usepackage{graphicx}
\usepackage{epsf}
\usepackage{pstricks}
\xyoption{all}
\usepackage{url}

\def\Box{\square}  

\def\mapright#1{\smash{\mathop{\longrightarrow}\limits^{#1}}}
\def\dmapright#1{\smash{\mathop{\Rightarrow}\limits^{#1}}}
\def\tra#1{\smash{\mathop{\mid\kern
-1pt\joinrel\relbar\joinrel\relbar}\limits^{*}_{#1}}}
\def\longtra#1{\smash{\mathop{\mid\kern
-1pt\joinrel\relbar\joinrel\relbar\joinrel\relbar}\limits^{*}_{#1}}}
\def\vlongtra#1{\smash{\mathop{\mid\kern
-1pt\joinrel\relbar\joinrel\relbar\joinrel\relbar\joinrel\relbar}\limits^{*}_{#1}}}
\def\vvlongtra#1{\smash{\mathop{\mid\kern
-1pt\joinrel\relbar\joinrel\relbar\joinrel\relbar\joinrel\relbar\joinrel\relbar}\limits^{*}_{#1}}}
\def\vvvlongtra#1{\smash{\mathop{\mid\kern
-1pt\joinrel\relbar\joinrel\relbar\joinrel\relbar\joinrel\relbar\joinrel\relbar\joinrel\relbar}\limits^{*}_{#1}}}
\def\etra#1{\smash{\mathop{\mid\kern
-1pt\joinrel\relbar\joinrel\relbar}\limits_{#1}}}
\def\mapleft#1{\smash{\mathop{\longleftarrow}\limits^{#1}}}

\def\vlongrightarrow{\relbar\joinrel\longrightarrow}
\def\vvlongrightarrow{\relbar\joinrel\vlongrightarrow}
\def\vvvlongrightarrow{\relbar\joinrel\vvlongrightarrow}

\def\longmapright#1{\smash{\mathop{\vlongrightarrow}\limits^{#1}}}
\def\vlongmapright#1{\smash{\mathop{\vvlongrightarrow}\limits^{#1}}}
\def\vvlongmapright#1{\smash{\mathop{\vvvlongrightarrow}\limits^{#1}}}

\def\A{{\cal{A}}}
\def\X{{\cal{X}}}
\def\AA{\tilde A}

\def\Rw{\Rightarrow}
\def\oo{\overline}

\def\B{{\cal{B}}}
\def\C{{\cal{C}}}

\def\N{\mathbb{N}}

\def\id{\textsf{id}}
\def\sing{\mbox{Sing}}

\def\rat{\mathop{\textsf{Rat}}}

\def\alm{\mathop{\textsf{AlmB}}}

\def\IS{\mbox{IS}}
\def\aut{\mathop{\textsf{Aut}}}

\def\max{\mbox{max}}
\def\min{\mbox{min}}

\def\G{{\cal{G}}}

\def\Z{\mathbb{Z}}

\def\p{\varphi}

\let\epsilon\varepsilon
\def\inv{^{-1}}

\def\beq{\begin{equation}}
\def\eeq{\end{equation}}
\def\bi{\begin{itemize}}
\def\ei{\end{itemize}}

\def\cc{\textsl{cc}}
\def\hc{\textsl{hc}}
\def\hcfp{\textsl{hcfp}}
\def\shcfp{\textsl{shcfp}}
\def\source{\textsl{source}}
\def\sink{\textsl{sink}}

\newtheorem{T}{Theorem}[section]
\newcommand{\bt}{\begin{T}}
\newcommand{\et}{\end{T}}
\newcommand{\ftd}{$\square$\end{T}}
\newtheorem{thm}[T]{Theorem}

\newtheorem{Proposition}[T]{Proposition}
\newcommand{\bp}{\begin{Proposition}}
\newcommand{\ep}{\end{Proposition}}
\newcommand{\fpd}{$\square$\end{Proposition}}

\newtheorem{Lemma}[T]{Lemma}
\newcommand{\bl}{\begin{Lemma}}
\newcommand{\el}{\end{Lemma}}
\newcommand{\fld}{$\square$\end{Lemma}}

\newtheorem{Corol}[T]{Corollary}
\newcommand{\bc}{\begin{Corol}}
\newcommand{\ec}{\end{Corol}}
\newcommand{\fcd}{$\square$\end{Corol}}

\newtheorem{Result}[T]{Result}
\newcommand{\br}{\begin{Result}}
\newcommand{\er}{\end{Result}}
\newcommand{\frd}{$\square$\end{Result}}

\newtheorem{Example}[T]{Example}
\newcommand{\be}{\begin{Example}}
\newcommand{\ee}{\end{Example}}

\newtheorem{Problem}[T]{Problem}
\newcommand{\bq}{\begin{Problem}}
\newcommand{\eq}{\end{Problem}}

\newtheorem{Fact}[T]{Fact}
\newcommand{\bfa}{\begin{Fact}}
\newcommand{\efa}{\end{Fact}}

\newtheorem{Remark}[T]{Remark}

\newcommand{\proof}
   {\par\medbreak\noindent{\bf Proof}.\enspace}

\newcommand{\qed}{
$\Box$
\par\bigbreak}

\newcounter{commentcounter}

\def\abstract#1{\par\bigskip
\begingroup\small
\baselineskip=12truept
\begin{center}ABSTRACT\end{center}
\par\medskip\par\noindent
\null\hfill\hbox{\vbox{\hsize=5truein\noindent#1}}
\hfill\null\par\endgroup\par}

\def\calE{\mathcal{E}}
\def\calL{\mathcal{L}}
\let\phi\varphi
\def\lefg{\le_{\textsf{fg}}}

\begin{document}

\title{Automorphic orbits in free groups: words versus subgroups
\thanks{The first author acknowledges support from Project ASA
(PTDC/MAT/65481/2006) and C.M.U.P., financed by F.C.T. (Portugal)
through the programmes POCTI and POSI, with national and European
Community structural funds.  This paper was prepared while the second
author was a visiting professor in the CSE Department, IIT Delhi. Both
authors acknowledge support from the 
ESF project \textsc{AutoMathA}.}}
\author{Pedro V. Silva, \small{\url{pvsilva@fc.up.pt}}\\
{\small Centro de Matem\'{a}tica, Universidade do
Porto}\thanks{%
Centro de Matem\'{a}tica, Faculdade de Ci\^{e}ncias, Universidade do
Porto, R. Campo Alegre 687, 4169-007 Porto, Portugal}
\and
Pascal Weil, \small{\url{pascal.weil@labri.fr}}\\
{\small LaBRI, Universit\'e de Bordeaux and CNRS}\thanks{%
LaBRI, Universit\'e Bordeaux-1, 351 cours de la Lib\'eration, 33405
Talence Cedex, France}}
\date{31th March 2010}
\maketitle

\abstract{
We show that the following problems are decidable in a rank 2 free
group $F_2$: does a given finitely generated subgroup $H$ contain
primitive elements?  and does $H$ meet the orbit of a given word $u$
under the action of $G$, the group of automorphisms of $F_2$?
Moreover, decidability subsists if we allow $H$ to be a rational
subset of $F_2$, or alternatively if we restrict $G$ to be a
rational subset of the set of invertible substitutions (a.k.a.
positive automorphisms).  In higher rank, the following weaker problem
is decidable: given a finitely generated subgroup $H$, a word $u$ and
an integer $k$, does $H$ contain the image of $u$ by some $k$-almost
bounded automorphism?  An automorphism is $k$-almost bounded if at
most one of the letters has an image of length greater than $k$.
}

\begin{center}\small
\end{center}

\newpage

Orbit problems in general concern the orbit of an element $u$ or a
subgroup $H$ of a group $F$, under the action of a subset $G$ of $\aut
F$.  Conjugacy problems are a special instance of such problems, where
$G$ consists of the inner automorphisms of $F$.  In this paper, we
restrict our attention to the case where $F$ is the free group $F_A$
with finite basis $A$.

In this context, orbit problems were maybe first considered by
Whitehead \cite{Whitehead}, who proved that membership in the
orbit of $u$ under the action of $\aut F_A$ is decidable.  The
analogous result regarding the orbit of a finitely generated subgroup
$H$ was established by Gersten \cite{Gersten}.  Much literature
has been devoted as well to the case where $G = \langle\p\rangle$ is a
cyclic subgroup of $\aut F_A$, e.g. Myasnikov and Shpilrain's work
\cite{MVW} on finite orbits of the form
$\langle\p\rangle\cdot u$ and Brinkmann's recent proof
\cite{Brinkmann} of the decidability of membership in
$\langle\p\rangle\cdot u$.

The orbit problems considered in this paper are of the following form:
given an element $u \in F_A$, a finitely generated subgroup $H$ of
$F_A$ and a subset $G$ of $\aut F_A$, does $H$ meet the orbit of $u$
under the action of $G$; that is: does $H$ contain $\p(u)$ for some
automorphism $\p\in G$?  A particular instance of this problem, when
$G = \aut F_A$, is the question whether $H$ contains a primitive
element, since the set of primitive elements of $F_A$ is the
automorphic orbit of each letter $a\in A$. The latter problem was
recently solved by Clifford and Goldstein in full generality \cite{CG}.
 These problems were posed
to the second author by O. Bogopolski, and they appear as Problem F39
in the list of open problems on \url{grouptheory.info}.

Our main results state that these problems are decidable in the rank 2
free group $F_2$, if $G = \aut F_2$ (Theorem~\ref{dorbit}) or if $G$
belongs to a certain family of rational subsets of $\aut F_2$, which
includes the rational subsets of invertible substitutions (a.k.a.
positive automorphisms, which map each letter to a positive word) or
of inverses of invertible substitutions, see Sections~\ref{sec:
rational autom} and~\ref{sec: invertible}.  For these rational values
of $G$, we also show the decidability of subgroup orbit problems:
given two finitely generated subgroups $H,K$ of $F_2$, does there
exist $\mu\in G$ such that $K$ is contained in (resp.  equal to)
$\mu(H)$.

In free groups with larger rank, we are only able to decide a weaker
problem.  Say that an automorphism $\p$ of $F_A$ is $k$-almost bounded
if $|\p(a)| > k$ for at most one letter $a\in A$.  We show that given
$k > 0$, $u\in F_A$ and $H$ a finitely generated subgroup of $F_A$,
one can decide whether there exists a $k$-almost bounded automorphism
$\mu$ such that $\mu(u) \in H$.

Some of our results hold also if we replace the subgroup $H$ by a
rational subset of $F_A$.

We use two main methods.  Some of our main results can be derived from
general results on the decidability of the solvability of equations
with rational constraints in free groups (Diekert, Guti\'errez and
Hagenah \cite{DGH}, building on Makanin's famous result
\cite{Makanin}). This is an interesting application of equations 
in free groups.

For other results, we give a direct combinatorial proof.  We use a
particular factorization of the automorphism group $\aut F_2$
(Theorem~\ref{deco}) and a detailed combinatorial analysis of the
effect of certain simple automorphisms on the graphical representation
of the subgroup $H$ (the representation by means of so-called
Stallings foldings \cite{Stallings,KM}, see
Section~\ref{sec automata}).  The set of these automorphisms is
$\Sigma = \{\p_{a,ba}, \p_{b\inv,a\inv}, \p_{b,a}\}$ ($\p_{x,y}$ maps
generator $a$ to $x$ and generator $b$ to $y$).

This combinatorial analysis leads to the definition of a (large but
finite) automaton whose vertices are finite automata associated with
the Stallings automata of the subgroups in the $\Sigma^*$-orbit of
$H$.  The construction of this automaton exploits the fact that a
certain combinatorial parameter of Stallings automata (which we call
the number of singularities) is preserved under the action of
automorphisms in $\Sigma$.  And it is the possibility of reading these
actions on this finite automaton which yields our decidability results
for the cases where $G$ is a rational subset of $\Sigma^*$.
Invertible substitutions form a particular rational submonoid of
$\Sigma^*$.

Interesting intermediary results state that the set of primitive
elements in $F_2$ is a context-sensitive language
(Proposition~\ref{csl}) and that if $|A| = m$ and $v_1,\ldots,
v_{m-1}\in F_A$, then the set of elements $x$ such that $v_1,\ldots,
v_{m-1}, x$ form a basis of $F_A$ is a constructible rational set
(Proposition~\ref{compbase}).

\section{Preliminaries}

\subsection{Free groups}

Let $A$ denote a finite alphabet.  The \textit{free monoid on} $A$,
written $A^*$, is the set of all finite sequences of elements of $A$
(including the empty sequence, written 1), under the operation of
concatenation.  We also write $A^+$ for the set of non-empty sequences
of elements of $A$.

Let $A\inv$ be a disjoint set of formal inverses of $A$ and let $\AA =
A \cup A\inv$.  The operation $u \mapsto u\inv$ is extended to $\AA^*$
as usual, by letting $(a\inv)\inv = a$ and $(ua)\inv = a\inv u\inv$
for all $a\in A$ and $u\in \AA^*$.

The \textit{free group} on $A$ is the quotient $F_A$ of $\AA^*$ by the
congruence generated by the pairs $(aa\inv,1)$, $a \in \AA$, and we
write $\pi\colon \AA^* \to F_A$ for the canonical projection.  A word
$u\in \AA^*$ is \textit{reduced} if it does not contain a factor
$aa\inv$ ($a\in \AA$) and we denote by $R_A$ the set of reduced words.
We also say that $u \in R_A$ is {\em cyclically reduced} if $uu$ is
reduced as well.  And we denote by $\cc(u)$ the \textit{cyclic core}
of $u$, that is, the unique word such that $u$ is of the form $u =
v\inv \cc(u) v$ in $\tilde A^*$.

We write $u\mapsto \bar u$ the \textit{reduction map}, where $\bar u$
is the (uniquely defined) word obtained from $u$ by iteratively
deleting factors of the form $aa\inv$ ($a\in \AA$) until none is left.
It is well-known that the reduction map is well defined, and that the
restriction $\pi\colon R_A \to F_A$ is a bijection.  To simplify
notation, if $g\in F_A$, we also write $\bar g$ for the reduced word
such that $\pi(\bar g) = g$, and we let the \textit{length} of $g$ be
$|g| = |\bar g|$.

Given $X \subseteq F_A$, we denote by $\langle X\rangle$ the subgroup
of $F_A$ generated by $X$.  We also let $\aut F_A$ denote the
automorphism group of $F_A$.  If $\p \in \aut F_A$ and no confusion
arises, we shall denote also by $\p$ the corresponding bijection of
$R_A$.

Given $B \subseteq F_A$, we say that $B$ is a {\em basis} of $F_A$ if
the homomorphism from $F_B$ to $F_A$ induced by the inclusion map $B \to
F_A$ is an isomorphism. Equivalently, $B$ is a basis of $F_A$ if and only if $B
= \p(A)$ for some $\p \in \aut F_A$. The \textit{primitive} elements 
of $F_A$ are those that sit in some basis of $F_A$. It follows that 
the set of primitive elements of $F_A$ is the orbit of each letter 
$a\in A$ under the action of $\aut F_A$.

In much of this paper, we shall be discussing the free group on 2
generators.  We fix the alphabet $A_2 = \{ a,b\}$ and use the notation
$F_2 = F_{A_2}$, $R_2 = R_{A_2}$.

\subsection{Automata and rational subsets}\label{sec automata}

The product of two subsets $K,L$ of a monoid $M$ is the subset $KL =
\{xy \mid x\in K,\ y\in L\}$.  The {\em star} operator on subsets is
defined by $L^* = \bigcup_{n\geq 0} L^n$, where $L^0 = \{1\}$.  A
subset $L$ of a monoid is said to be {\em rational} if $L$ can be
obtained from finite subsets using finitely many times the operators
union, product and star.  We denote by $\rat M$ the set of rational
subsets of $M$.  If $M$ is the free monoid $A^*$ on a finite alphabet
$A$, subsets of $A^*$ are called \textit{languages}, and elements of
$\rat A^*$ are called \textit{rational languages}.

Note that if $\p\colon A^* \to M$ is an onto morphism, then $\rat M$
is the set of all $\p(L)$ where $L\in \rat A^*$.  For instance, every
finitely generated subgroup of $F_A$ is rational.

It is well-known that rational languages can be characterized by means
of finite automata.  A (finite) $A$-automaton is a tuple $\A =
(Q,q_0,T,E)$ where $Q$ is a (finite) set, $q_0 \in Q$, $T \subseteq Q$
and $E \subseteq Q \times A \times Q$.  It can be viewed as a graph
with vertex set $Q$ (the \textit{states}), with a designated vertex
$q_0$ (the \textit{initial state}) and a set of designated vertices
$T$ (the \textit{terminal} states), whose edges are labeled by letters
in $A$, and are given by the set $E$ (the \textit{transitions}).

A {\em nontrivial path} in $\A$ is a sequence
$$p_0 \mapright{a_1} p_1 \mapright{a_2} \ldots \mapright{a_n} p_n$$
with $n\ge 1$, $(p_{i-1},a_i,p_i) \in E$ for $i = 1,\ldots,n$.  Its
{\em label} is the word $a_1\ldots a_n \in A^+$.  We consider also the
{\em trivial path} $p_0 \mapright{1} p_0$ for each $p_0 \in Q$, whose label
is the empty word.  A path is said to be {\em successful} if $p_0 =
q_0$ and $p_n \in T$.  The {\em language} $L(\A)$ {\em recognized by}
$\A$ is the set of all labels of successful paths in $\A$.

The automaton $\A = (Q,q_0,T,E)$ is said to be {\em deterministic} if,
for all $p \in Q$ and $a \in A$, there is at most one edge of the form
$(p,a,q)$.  In that case, we write $q = p\cdot a$.  We say that $\A$
is {\em trim} if every $q \in Q$ lies in some successful path.

Kleene's theorem states that a language is rational if and only if it
is accepted by a finite automaton, which can be required to be
deterministic and trim, see \cite{HU}.
In the context of a particular result or claim, we say that a rational
language $L$ is {\em effectively constructible} if there exists an
algorithm to produce a finite automaton recognizing $L$ from the
concrete structures containing the input.  More generally, if
$\p\colon A^* \to M$ is an onto morphism, we say that a rational
subset of $M$ is \textit{effectively constructible} (with respect to
$A$) if it is the image of an effectively constructible rational language
over $A$.

\begin{Remark}\rm
    A subset $L\subseteq F_A$ is rational if $L =
    \pi(K)$ for some rational subset $K$ of $\tilde A^*$.  Benois'
    theorem \cite{Benois} states that this is the case if and only
    if $\oo{K}$ ($= \pi\inv(\pi(K)) \cap R_A$, in bijection with $L$
    via $\pi$) is a rational subset of $\tilde A^*$. In the sequel we 
    sometimes confuse the notions of a rational subset of $F_A$ and a 
    rational language in $\tilde A^*$ that consists only of reduced 
    words.
\end{Remark}

\subsection{Automata and subgroups of $F_A$}

To discuss subgroups of free groups, we use inverse automata.  In an
$\AA$-automaton $\A = (Q,q_0,T,E)$, the {\em dual} of an edge
$(p,a,q) \in E$ is $(q,a\inv,p)$.  Then $\A$ is said to be {\em dual} if $E$
contains the duals of all edges, and {\em inverse} if it is dual,
deterministic, trim (equivalent to {\em connected} in this case) and $|T| = 1$.

Given a finitely generated subgroup $H$ of $F_A$ (we write $H \lefg
F_A$), we denote by $\A(H)$ the \textit{Stallings automaton}
associated to $H$ by the construction often referred to as {\em
Stallings foldings}.  This construction, that can be traced back to
the early part of the twentieth century \cite[Chapter 11]{Rotman},
was made explicit by Serre \cite{Serre} and Stallings
\cite{Stallings} (see also \cite{KM}).

A brief description is as follows.  If $h_1,\ldots,h_r \in R_A$ is a
set of generators of the subgroup $H$ (that is, $H = \langle
\pi(h_1),\ldots,\pi(h_r)\rangle$), one constructs a dual automaton in
the form of $r$ subdivided circles around a common distinguished
vertex 1 (both initial and terminal), each labeled by one of the
$h_i$.  Then we iteratively 
identify identically labeled pairs of edges starting (resp.  ending)
at the same vertex (this is called the \textit{folding} process, see
Figure~\ref{figure folding}), until no further folding is possible.

\begin{figure}
    $$\xymatrix{
&{q} &&&\\
{p} \ar^a[ur] \ar_a[r]& {r} && {p}\ar_a[r] & {q\sim r}
}$$
\caption{A folding step, with $a\in \AA$}\label{figure folding}
\end{figure}

The following proposition summarizes important properties, see
\cite{KM}.

\bp
\label{RSS}
Let $H \lefg F_A$.  Then:
\begin{itemize}
    \item[(i)] $\A(H)$ is a finite inverse automaton, which does not
    depend on the finite reduced generating set nor on the 
    sequence of foldings chosen;
    
    \item[(ii)] if $p\mapright{u} q$ is a path in $\A(H)$, so is $p
    \mapright{\oo{u}} q$;
    
    \item[(iii)] for every $u \in R_A$, $u \in L(\A(H))$ if and only
    if $\pi(u) \in H$; in particular, $L(\A(H)) \subseteq \pi\inv(H)$;

    \item[(iv)] for every cyclically reduced $u \in F_A$, $wuw\inv \in
    H$ for some $w \in F_A$ if and only if $u$ labels some loop in
    $\A(H)$.
\end{itemize}
\ep

\section{The mixed orbit problem as an equation problem}
\label{sec: orbit pbms equation}

Our original motivation on writing this paper was solving the mixed
orbit problem $\p(u) \in H$ for given $u \in F$ and $H \lefg F$.
We managed to solve it in rank 2, see Section~\ref{sec: equations}
below.  We also solve it in arbitrary rank, if we impose a restriction
on the class of automorphisms, see Section~\ref{beyond rank 2} below.

Our initial proof in rank 2 was purely combinatorial, and is presented
in Section~\ref{sec: proof} below, whereas our result in higher
ranks made use of a result of Diekert, Guti\'errez and Hagenah
\cite{DGH} on equations with rational constraints which we will
discuss below.  Upon reading a first version of this
paper\footnote{arXiv:0809.4386v1 [math.GR]}, Dahmani and Girardel, and
independently Enric Ventura (recalling a foregone conversation with
Alexei Miasnikov), called our attention to the fact that Diekert,
Guti\'errez and Hagenah's result could also be used to prove our
result in rank 2.

\subsection{Equations with rational constraints and automorphisms of $F_2$}
\label{sec: equations}

A \textit{system of equations} in a free group $F_A$, with set of
unknowns $X$ (disjoint from $A$), is a tuple $\calE =
(e_1,\ldots,e_k)$ of elements of $F_{A\cup X}$.  A \textit{solution}
of $\calE$ is a morphism from $F_{A\cup X}$ to $F_A$ which maps each
letter of $A$ to itself and every element of $\calE$ to $1$.

A \textit{rational constraint} on a system of equations with unknowns
in $X$ is a collection $\calL = (L_x)_{x\in X}$ of rational subsets of
$F_A$ and we say that a morphism $\phi\colon F_{A\cup X} \to F_A$ is a
solution of the system $\calE$ with rational constraints $\calL$ if
$\phi$ is a solution of $\calE$ and $\phi(x) \in L_x$ for each $x\in
X$.  Diekert, Guti\'errez and Hagenah \cite{DGH} showed the following
result.

\begin{thm}\label{thm: DGH}
    The satisfiability problem for systems of equations with rational 
    constraints in a free group is decidable.
\end{thm}

This leads to a quick solution of the mixed orbit problem in $F_2$.

\begin{thm}\label{mixed orbit rational}
    Given $u \in F_2$ and a rational subset $L$ of $F_2$, it is
    decidable whether or not $\p(u) \in L$ for some $\p \in \aut
    F_2$.
\end{thm}

\proof
By a result attributed to Dehn, Magnus and Nielsen (see
\cite{Shpilrain}), $\{x,y\}$ is a basis of $F_2$ if and only if
there exists some $g \in F_2$ such that $g\inv[x,y]g = [a,b]^{\pm 1}.$

Now observe that if $u\in F_2$, then $\p(u) \in H$ for some $\p \in
\aut F_2$ if and only if $u(x,y) \in H$ for some basis $\{x,y\}$ ---
where $u(x,y)$ denotes the word $u$ in which each occurrence of $a$ 
has been replaced by $x$ and each occurrence of $b$ by $y$.

Thus there exists an automorphism $\p$ such that $\p(u) \in H$ if and 
only if one of the following systems (in
the unknowns $x,y,v,g$, and for $\epsilon = \pm 1$) admits a solution 
$$\left\{
\begin{array}{l}
g\inv[x,y]g = [a,b]^{\epsilon}\\
u(x,y) = v
\end{array}
\right.$$
with the rational constraint that $v\in L$.  This is decidable by
Theorem~\ref{thm: DGH}.
\qed

Since a finitely generated subgroup is rational (it is the star of 
its generators and their inverses), Theorem~\ref{mixed orbit 
rational} yields the following corollary.

\begin{Corol}\label{dorbit}
    Given $u \in F_2$ and $H \lefg F_2$, it is decidable whether or
    not $\p(u) \in H$ for some $\p \in \aut F_2$.
\end{Corol}

Since the primitive elements of $F_2$ are the orbit of each letter 
$a\in A$ under $\aut F_2$, we also note the following result.

\begin{Corol}\label{primitive}
    Given a rational subset $L$ of $F_2$ (e.g. a finitely generated
    subgroup), it is decidable whether or not $L$ contains a primitive
    element.
\end{Corol}

\begin{Remark}\rm
    Clifford and Goldstein also proved a comparable result for primitive 
    elements, by completely different methods: they show that it is 
    decidable whether a finitely generated subgroup $H$ of $F_A$ (for 
    any finite alphabet $A$) contains a primitive element 
    \cite{CG}.
\end{Remark}

One can also consider, in the statement of Theorem~\ref{mixed orbit
rational}, the rational subset of positive elements of $F_2$ (namely,
the submonoid $A^*$).  Then our result shows that it is decidable
whether an element $u \in F_2$ is \textit{potentially positive}, that
is, whether it has a positive automorphic image.  Different proofs of
this result already appear in Goldstein \cite{Goldstein} and Lee
\cite{Lee}.

Another idea is to consider a tuple of elements of $F_2$ rather than a
single element $u$, or equivalently a subgroup $K\lefg F_2$.

\begin{thm}\label{extensions equations F2}
    Let $u_1,\ldots,u_k\in F_2$, $L_1,\ldots,L_k\in\rat F_2$ and $H,K \lefg F_2$.
    The following problems
    are decidable:
    \begin{itemize}
	\item[(1)] whether $\p(u_1) \in L_1,\ldots,\p(u_k) \in L_k$, for
	some $\p\in\aut F_2$;
	
	\item[(2)] whether conjugates of $\p(u_1),\ldots,\p(u_k)$ sit
	in $L_1, \ldots, L_k$, respectively, for some $\p\in\aut F_2$;
	
	\item[(3)] whether $\p(K) \subseteq H$, for some
	$\p\in\aut F_2$.
    \end{itemize}
\end{thm}

\proof
These statements are proved like Theorem~\ref{mixed orbit rational}, 
by reduction to a system of equations with rational constraints.

For (1), we consider the system
$$\left\{
\begin{array}{l}
g\inv[x,y]g = [a,b]^{\epsilon}\\
u_1(x,y) = v_1\\
\hskip.5cm\ldots\\
u_k(x,y) = v_k
\end{array}
\right.$$
in the unknowns $x,y,g,v_1,\ldots,v_k$ with the rational constraints
$v_1\in L_1,\ldots,v_k\in L_k$.

For (2), we consider the system
$$\left\{
\begin{array}{l}
g\inv[x,y]g = [a,b]^{\epsilon}\\
h_1\inv u_1(x,y)h_1 = v_1\\
\hskip.5cm\ldots\\
h_k\inv u_k(x,y)h_k = v_k
\end{array}
\right.$$
in the unknowns $x,y,g,h_1,\ldots,h_k,v_1,\ldots,v_k$ with the
rational constraints 
$v_1\in L_1,\ldots,v_k\in L_k$.

Statement (3) is a particular case of (1), when $u_1,\ldots,u_k$ is a 
generating set of $K$ and $L_1 = \ldots = L_k  = H$.
\qed

\begin{Remark}\label{remark complexity}
    \rm Instead of asking whether there exists an automorphism in
    $\aut F_2$ mapping $u$ into $L$, one may want to exhibit such an
    automorphism, if it exists.
    
    The existence question reduces to the satisfiability of equations
    with rational constraints, and Diekert, Guti\'errez and Hagenah
    showed that this can be done in PSPACE \cite{DGH}.  One
    can extract from that paper a description of such a solution (an
    automorphism) as an exponential length product of simple
    automorphisms: the images of the letters may therefore have double
    exponential length.
\end{Remark}

\subsection{Beyond rank 2} \label{beyond rank 2}

We do not know how to extend Theorem~\ref{mixed orbit rational} or
Corollary \ref{dorbit} to arbitrary finite alphabets, but we can get
decidability for weakened versions of the problem.  The first such
result involves a restriction on the subgroups considered.

\bt
Let $u \in F_A$ and let $H \lefg F_A$.  If $H$ is cyclic or a
free factor of $F_A$, it is decidable whether or not $\p(u) \in H$
for some $\p \in \aut F_A$.
\et

\proof
Let us first assume that $H$ is a free factor of $F_A$, with rank $k$.
It is easily verified that $\p(u) \in H$ for some automorphism $\p$ if
and only if $u$ sits in some rank $k$ free factor of $F_A$.  This is
known to be decidable: it suffices to compute a minimum length element
$v$ in the automorphic orbit of $u$ (the so-called easy part of
Whitehead's algorithm, see \cite{LS,RVW}) and to
verify whether $v$ uses at least $k$ letters (Shenitzer
\cite{Shenitzer}, see also \cite[Prop.  I.5.4]{LS}).

Let us now assume that $H = \langle v\rangle$.  Without loss of
generality, we may assume that $u$ and $v$ are cyclically reduced.
Say that a word $x$ is root-free if it is not equal to a non-trivial
power of a shorter word.  Then $u = x^k$ for some uniquely determined
integer $k\ge 1$ and root-free word $x$, and similarly, $v = y^\ell$
for some uniquely determined $\ell \ge 1$ and root-free $y$.  It is an
elementary verification that the image of a 
root-free word by an automorphism is also 
root-free.  Thus, an automorphism maps $u$ into $H$ if and only if 
 it maps $x$ to $y$ or $y\inv$, and $k$ is a multiple of $\ell$.
Decidability follows from the fact that we can decide whether two
given words are in each other's automorphic orbit, using Whitehead's
algorithm \cite{LS}.
\qed

The second result on a weakened version of our orbit problem involves
{\em almost bounded automorphisms}.  Given a finite alphabet $A$ and
$k \in \N$, we say that an automorphism $\p$ of $F_A$ is
\textit{$k$-almost bounded} if $|\p(a)| > k$ for at most one letter $a
\in A$.  We let $\alm_k F_A$ denote the set of $k$-almost bounded
automorphisms of $F_A$.

\bt
\label{almost}
Given $u \in F_A$, $L \in \rat F_A$ and $k \in \N$, it is decidable
whether or not $\p(u) \in L$ for some $\p \in \alm_k F_A$.
\et

The proof of this theorem relies on Diekert, Guti\'errez and Hagenah's
result on the satisfiability of equations with rational constraints in
free groups discussed in Section~\ref{sec: equations}.  We also
require two technical results.

\begin{Lemma}\label{basa}
    Let $A = \{ a_1,\ldots,a_m\}$ and $u \in R_A$.  Then $\{
    a_1,\ldots,a_{m-1},u\}$ is a basis of $F_A$ if and only if $u =
    va_m^{\epsilon}w$ for some $v,w \in R_{\{a_1,\ldots,a_{m-1}\}}$
    and $\epsilon \in \{ 1,-1\}$.
\end{Lemma}

\proof
It is immediate that if $u = v a_m^\epsilon w$ with $v, w \in
R_{\{a_1,\ldots,a_{m-1}\}}$, then $\{a_1,...,a_{m-1},u\}$ generates
$F_A$, and by the Hopfian property of free groups (see \cite[Prop.
I.3.5]{LS}), $\{a_1,...,a_{m-1},u\}$ is a basis of $F_A$.

Conversely, let $u \in R_A$ contain at least an occurrence of $a_m$ or
$a_m\inv$, and let $u = vzw$ be the factorization with $v, w \in
R_{\{a_1,...,a_{m-1}\}}$ of maximal length.  It is immediate that if
$H = \langle a_1,...,a_{m-1},u\rangle$, then $H = \langle
a_1,...,a_{m-1},z\rangle$ and $\A(H)$ is equal to $\A(\langle
z\rangle)$ with loops labelled $a_1,...,a_{m-1}$ attached at the
origin.  Thus, if $\{a_1,...,a_{m-1},u\}$ is a basis of $F_A$, then
$\A(\langle z\rangle)$ must consist of a single loop labeled $a_m$,
and hence $z$ must be equal to $a_m$ or $a_m\inv$.
\qed

This leads to the following generalization.

\bp
\label{compbase}
Let $m = |A|$ and $v_1,\ldots,v_{m-1} \in R_A$. Then
$$X = \{ x \in R_A\mid (v_1,\ldots,v_{m-1},x)\mbox{ is a basis of
}F_A\}$$
is rational and effectively constructible.
\ep

\proof
First note that $X$ is nonempty if and only if $(v_1,...,v_{m-1})$ is
a basis of a free factor of $F_A$.  This is well-known to be
decidable.  Moreover, if $X \ne\emptyset$, then we can effectively
construct an element $z$ of $X$: it is verified in
\cite{SW} that if $K = \langle v_1,..., v_{m-1}\rangle$,
then $K$ is a free factor of $F_A$ if and only if there are vertices
$p$ and $q$ of $\A(K)$ whose identification leads (via foldings) to
the bouquet of circles $\A(F_A)$, and in that case, if $u_p$ and $u_q$
are the labels of geodesic paths of $\A(K)$ from the origin to $p$ and
$q$, then $z = \overline{u_p u_q^{-1}} \in X$.

Let $\p \in \aut F_A$ be defined by $\p(a_i) = v_i$ ($i = 1,\ldots,
m-1$) and $\p(a_m) = z$.  Then $x \in X$ if and only if
$(a_1,\ldots,a_{m-1},\p\inv(x))$ is a basis of $F_A$.  By Lemma
\ref{basa}, this is equivalent to say that $\p\inv(x) \in R(a_m \cup
a_m\inv)R$, where $R = \langle a_1,\ldots,a_{m-1}\rangle$, and
therefore
$$X = \p(R(a_m \cup a_m\inv)R) = V(z \cup z\inv)V$$
for $V = \langle v_1,\ldots,v_{m-1}\rangle$.

In particular, $X$ is rational and the formula $X = V(z \cup z\inv) V$
provides an effective construction for it.
\qed

\par\medbreak\noindent{\bf  Proof of Theorem \ref{almost}}. 
Write $A = \{ a_1,\ldots,a_m\}$.  Without loss of generality, we may
restrict ourselves to the case $|\p(a_i)|\; \leq k$ for $i =
1,\ldots,m-1$.  Since there are only finitely many choices for these
$\p(a_i)$, we may as well assume them to be fixed, say $\p(a_i) =
v_i$ for $i = 1,\ldots,m-1$.  Let then $X = \{ x \in R_A\mid
(v_1,\ldots,v_{m-1},x)\mbox{ is a basis of }F_A\}$: then $X$ is
rational by Proposition \ref{compbase}.

Write $u = u_0a_m^{\epsilon_1}u_1 \ldots a_m^{\epsilon_n}u_n$ with $n
\geq 0$, $u_i \in F_{\{a_1,\ldots,a_{m-1}\}}$ and $\epsilon_i = \pm 1$
for every $i$.  Then we must decide whether there exists some $y \in
X$ such that
$$u'_0y^{\epsilon_1}u'_1 \ldots y^{\epsilon_n}u'_n \in L,$$
where $u'_i = u_i(v_1,\ldots,v_{m-1})$ is the word obtained from $u_i$
by replacing each $a_j$ by $v_j$.  This is equivalent to deciding
whether the equation
\begin{equation}\label{eq almost bounded}
    u'_0y^{\epsilon_1}u'_1 \ldots y^{\epsilon_n}u'_n = z
\end{equation}
on the variables $y,z$ has a solution in $F_A$ with the rational
constraints $y \in X$ and $z \in L$.  This is decidable by
Theorem~\ref{thm: DGH}.
\qed

As in Section~\ref{sec: equations}, Theorem \ref{almost} yields a 
decidability result for finitely generated subgroups.

\bc
\label{coralmost}
Given $u \in F_A$, $H \lefg F_A$ and $k \in \N$, it is decidable
whether or not $\p(u) \in H$ for some $\p \in \alm_k F_A$.
\ec

And as in Theorem~\ref{extensions equations F2}, we use the same ideas
to prove the following theorem.  If $w\in F_2$, $\lambda_w$ denotes
the inner automorphism $u \mapsto w\inv uw$.

\begin{thm}
    Let $u_1,\ldots,u_m\in F_A$, $L\in\rat F_A$, $H,K \lefg F_A$ and
    $k\in\N$.  The following problems are decidable:
    \begin{itemize}
	\item[(1)] whether $\p(u_1),\ldots,\p(u_m) \in L$, for
	some $\p\in\alm_k F_A$;
	
	\item[(2)] whether $\lambda_w\p(u_1),\ldots,\lambda_w\p(u_m)
	\in L$, for some $w\in F_A$ and $\p\in\alm_k F_A$;
	
	\item[(3)] whether conjugates of $\p(u_1),\ldots,\p(u_m)$ sit
	in $L$, for some $\p\in\alm_k F_A$;
	
	\item[(4)] whether $\p(K) \subseteq H$, for some
	$\p\in\alm_k F_A$;
	
	\item[(5)] whether $\lambda_w\p(K) \subseteq H$, for some
	$w\in F_A$ and $\p\in\alm_k F_A$.
    \end{itemize}
\end{thm}

\proof
These statements are proved like Theorem~\ref{almost}, by reduction to
a system of equations with rational constraints.

For the first statement, we consider a system of equations of the form
of equation~(\ref{eq almost bounded}) in the proof of
Theorem~\ref{almost}, one for which $u_j$, $1\le j\le m$ (the unknowns
are $y, z_1,\ldots,z_m$).

For the second (resp. third) statement, we consider the same system, with each 
equation conjugated by a new unknown $v$ (resp. by distinct new unknowns 
$v_j$, $1\le j\le m$).

The fourth and fifth statements are applications of the first and
second when the rational subset $L$ is the subgroup $H$.
\qed

\begin{Remark}
    \rm Following-up with the discussion in Remark~\ref{remark
    complexity}, we note that the complexity upper bounds for the
    decision problems described in this section are PSPACE again: we
    need to (attempt to) solve, successively, systems of equations for
    the different values of $v_1,\ldots,v_{m-1}$ (with the notation of
    the proof of Theorem~\ref{almost}.  These $(m-1)$-tuples of words
    of length at most $k$ are exponentially many (in the variable $k$)
    but they can be listed in polynomial space.
\end{Remark}

\section{Combinatorial approach: the role of $\Sigma$}

We now restrict our attention to $F_2$. 
If $x,y\in F_2$, we denote by 
$\p_{x,y}$ the endomorphism mapping $a$ to $x$ and $b$ to $y$. 
In this section, we discuss some properties of the following sets of 
automorphisms of $F_2$:
\begin{itemize}
    \item[] $\Sigma_0 = \big\{ \p_{a,ba}, \p_{b\inv,a\inv} \big\}$ and
    $\Sigma = \Sigma_0 \cup \big\{ \p_{b,a} \big\}$;
    
    \item[] $\Phi = \big\{ \p_{a,ba}, \p_{ab,b},\p_{a,ab}, \p_{ba,b}
    \big\}$;

    \item[] $\Delta = \big\{ \p_{a,a^mb^{\epsilon}a^n} \mid \; m,n \in \Z,\;
    \epsilon \in \{ 1,-1\} \big\}$;
    
    \item[] $\Psi =\big\{ \p \in \aut F_2 : |\p(a)| = |\p(b)| = 1
    \big\}$ and $\Lambda = \big\{ \lambda_w \mid w \in R_A\big\}$.
\end{itemize}

The following will be useful in the sequel.

\bp 
\label{proptec}
\bi
\item[(i)]
$X\Lambda = \Lambda X$ for every $X \subseteq \aut F_2$;
\item[(ii)]
$\Lambda\Psi\Phi^* \subseteq \Lambda\Psi(\Sigma_0\inv)^*\p_{a\inv,b}$;
\item[(iii)]
$\Delta \subseteq \Lambda(\p_{a,ba}^*\cup
\p_{a\inv,b}\p_{a,ba}^*\p_{a\inv,b})(1 \cup \p_{a,b\inv})$. 
\ei
\ep

\proof
(i) follows from the fact that $\theta\lambda_w =
\lambda_{\theta(w)}\theta$ for each $w\in F_2$ and $\theta\in\aut
F_2$.

(ii) Notice that $\p_{ab,b} = \p_{b,a}\p_{a,ba}\p_{b,a}$, $\p_{a,ab} = 
\lambda_{a\inv}\p_{a,ba}$ and $\p_{ba,b} = \lambda_{b\inv}\p_{ab,b}$. It 
follows that $\Lambda\Psi\Phi^* \subseteq \Lambda\Psi \{\p_{a,ba}, 
\p_{b,a}\}^*$.

Observe also that $\p_{a,ba} = \p_{a\inv,b}\p_{a,ba}\inv\p_{a\inv,b}$,
$\p_{b,a} = \p_{a\inv,b}\p_{b\inv,a\inv}\inv\p_{a\inv,b}$ and
$\p_{a\inv,b}^2 = 1$.  So we have
$$\{\p_{a,ba}, \p_{b,a}\}^* =
\p_{a\inv,b}\{\p_{a,ba}\inv,\p_{b\inv,a\inv}\inv\}^*\p_{a\inv,b} =
\p_{a\inv,b}(\Sigma_0\inv)^*\p_{a\inv,b}.$$
Therefore $\Lambda\Psi\Phi^* \subseteq
\Lambda\Psi\p_{a\inv,b}(\Sigma_0\inv)^*\p_{a\inv,b} =
\Lambda\Psi(\Sigma_0\inv)^*\p_{a\inv,b}$.

(iii) Observe that if $m,n\in \Z$, then $\p_{a,a^mba^n} =
\lambda_{a^{-m}}\p_{a,ba^{m+n}} = \lambda_{a^{-m}}\p_{a,ba}^{m+n}$, so that
$\p_{a,a^mba^n} \in \Lambda (\p_{a,ba}^* \cup (\p_{a,ba}\inv)^*)$.  We
already noted that $\p_{a,ba}\inv = \p_{a\inv,b} \p_{a,ba}
\p_{a\inv,b}$ and $\p_{a\inv,b}^2 = 1$, so
$$\p_{a,a^mba^n} \in \Lambda(\p_{a,ba}^* \cup
\p_{a\inv,b}\p_{a,ba}^*\p_{a\inv,b}).$$
Similarly, $\p_{a,a^mb\inv a^n} =
\lambda_{a^{n}}\p_{a,ba}^{-(m+n)}\p_{a,b\inv}$ and hence
$$\p_{a,a^mb\inv a^n} \in \Lambda(\p_{a,ba}^* \cup
\p_{a\inv,b}\p_{a,ba}^*\p_{a\inv,b})\p_{a,b\inv},$$
which concludes the proof.
\qed

\subsection{Primitive words and a factorization of $\aut F_2$}

Let us first consider a particular automorphic orbit in $F_A$, namely
the set $P_A$ of primitive words.  Recall that a word is {\em
primitive} if it belongs to some basis of $F_A$.  In particular, $P_A$
is the automorphic orbit of each letter from $A$.  We shall often view
$P_A$ as a subset of $R_A$.  We denote by $P_2$ the set of all
primitive words in $F_2$.

We use a known characterization of the words in $P_2$ to derive a
technical factorization of the group $\aut F_2$ of automorphisms of
$F_2$, that will be used in Section~\ref{sec: orbit pbms}.  We further
exploit this characterization to point out certain language-theoretic
properties of $P_2$.

Proposition~\ref{new34} reports two results: the first is due to
Nielsen \cite{Nielsen} (see also \cite[2.2]{Cohen} and \cite{OZ81})
and the second is due to Wen and Wen \cite{WenWen}.  An interesting
perspective on either is offered in \cite[Chapter 2]{Lothaire2} and
\cite[Chapter I-5]{Christoffel}.

\bp
\label{new34}
\bi
\item[(i)] Up to conjugation, every primitive element $u
\in P_2$ is either a  
letter, or of the form $u = a^{n_1}b^{m_1} ... a^{n_k}b^{m_k}$ where
\bi
\item[-]
either $n_1 = ... = n_k \in \{1,-1\}$ and $\{m_1,...,m_k\} \subseteq
\{n,n+1\}$ for some integer $n$,
\item[-]
or $m_1 = ...  = m_k \in \{1,-1\}$ and $\{n_1,...,n_k\} \subseteq
\{n,n+1\}$ for some integer $n$.
\ei
\item[(ii)]
The set of positive primitive words $P_2 \cap \{a,b\}^+$ is equal to
$\Phi^*(\{a,b\}) = b \cup \Phi^*(a)$.
\ei
\ep

\bc
\label{new35}
$P_2 = \Lambda\Psi\Phi^*(a)$.
\ec

\proof
By Proposition \ref{new34} (i), every primitive element of $F_2$ is a
conjugate of $\psi(ab^{m_1} ...  ab^{m_k})$, where $\{m_1,...,m_k\}
\subseteq \{n,n+1\}$ for some integer $n\ge 0$ and $\psi \in \Psi$.
That is, $P_2 = \Lambda\Psi(P_2 \cap \{ a,b\}^+)$.  By Proposition
\ref{new34} (ii), it follows that $P_2 = \Lambda\Psi(b \cup \Phi^*(a))
= \Lambda\Psi\Phi^*(a)$.
\qed

We can now prove a useful decomposition result for $\aut F_2$.

\bt
\label{deco}
$\aut F_2 = \Lambda\Psi\Phi^*\Delta = \Psi(\Sigma_0\inv)^*\Lambda
\p_{a,ba}^*(\p_{a\inv,b}\cup \p_{a\inv,b\inv})$.
\et

\proof
To establish the first equality, we consider $\theta \in \aut F_2$.
Then $\theta(a) \in P_2$ and so $\theta(a) = \sigma(a)$ for some
$\sigma \in \Lambda\Psi\Phi^*$ by Corollary \ref{new35}.  Corollary
\ref{basa} then shows that $\sigma\inv\theta =
\p_{a,a^mb^{\epsilon}a^n}$ for some $m,n \in \Z$ and $\epsilon \in \{
1,-1\}$.  So $\sigma\inv\theta \in \Delta$ and $\theta \in
\Lambda\Psi\Phi^*\Delta$. It follows that
\begin{align*}
    \aut F_2 &\subseteq \Lambda\Psi\Phi^*(\p_{a,ba}^*\cup
    \p_{a\inv,b}\p_{a,ba}^*\p_{a\inv,b})(1 \cup \p_{a,b\inv})\textrm{ by 
    Proposition \ref{proptec}}\\
    & \subseteq \Lambda\Psi\Phi^*(1\cup
    \p_{a\inv,b}\p_{a,ba}^*\p_{a\inv,b})(1 \cup \p_{a,b\inv})\textrm{ 
    since $\p_{a,ba} \in \Phi$}\\
    & \subseteq
    \Lambda\Psi\Phi^*(\p_{a\inv,b}\p_{a,ba}^*\p_{a\inv,b})(1 \cup
    \p_{a,b\inv})\textrm{ since $\p_{a\inv,b}^2 = 1$}\\
    &\subseteq \Lambda\Psi(\Sigma_0\inv)^* \p_{a,ba}^* \p_{a\inv,b} (1
    \cup \p_{a,b\inv})\textrm{ by Proposition \ref{proptec} (ii)}\\
    & \subseteq \Psi(\Sigma_0\inv)^*\Lambda
    \p_{a,ba}^*(\p_{a\inv,b}\cup \p_{a\inv,b\inv}).
\end{align*}
The converse inclusion is of course trivial.
\qed

\subsection{Invertible substitutions}\label{sec: invertible}

A {\em substitution}\footnote{also called a \textit{positive
endomorphism}} of $F_A$ is an endomorphism $\p$ such that $\p(a) \in
A^*$ for every $a \in A$.  If $\p$ is an automorphism, it is said to
be an {\em invertible substitution}.  We denote by $\IS(F_2)$ the
monoid of all invertible substitutions of $F_2$, and by $\IS\inv(F_2)$
the monoid of their inverses. Note that the 
inverse of an invertible substitution is not necessarily a 
substitution: indeed $\p_{a,ba}\inv = \p_{a,ba\inv}$.

\begin{Lemma}\label{lemma: IS and ISinv}
    $\IS(F_2)$ is a rational submonoid of $\Sigma^*$. Moreover, there 
    exists a rational submonoid $S$ of $\Sigma^*$ such that
    $\IS\inv(F_2) = \p_{a,b\inv}\ S\ \p_{a,b\inv}$.
    
    In addition, every rational subset $R \in \rat \IS(F_2)$ is also 
    in $\rat \Sigma^*$, and every rational subset $R \in \rat 
    \IS\inv(F_2)$ is of the form $\p_{a,b\inv}\ R'\ \p_{a,b\inv}$ for 
    some $R'\in \rat S \subseteq \rat\Sigma^*$.
\end{Lemma}
    
\proof
It is known \cite{WenWen} that the monoid $\IS(F_2)$ is generated by
$\p_{b,a}$, $\p_{a,ba}$ and $\p_{a,ab}$ (see also \cite[Chapter
I.5]{Christoffel}, \cite[Sec.  2.3.5]{Lothaire2}).  But $\p_{b,a},
\p_{a,ba} \in \Sigma$ and
$$\p_{a,ab} = 
\p_{b,a}\p_{b\inv,a\inv}\p_{a,ba}\p_{b\inv,a\inv}\p_{b,a} \in 
\Sigma^*,$$
so $\IS(F_2) = \{\p_{b,a}, \p_{a,ba}, \p_{a,ab}\}^* \in \rat
\Sigma^*$.

Next we observe that
\begin{align*}
    \p_{b,a}\inv &= \p_{b,a} =
    \p_{a,b\inv}\p_{b\inv,a\inv}\p_{a,b\inv} \\
    \p_{a,ba}\inv &= \p_{a,ba\inv} = 
    \p_{a,b\inv}\p_{a,ab}\p_{a,b\inv}\\
    \p_{a,ab}\inv &= \p_{a,a\inv b} = 
    \p_{a,b\inv}\p_{a,ba}\p_{a,b\inv}
\end{align*}
Since $\p_{a,b\inv}$ has order 2, it follows that $\IS(F_2)\inv =
\p_{a,b\inv}R\p_{a,b\inv}$ with $R = \{\p_{b\inv,a\inv}, \p_{a,ab},
\p_{a,ba}\}^* \in \rat \Sigma^*$.
\qed

\subsection{Primitive words form a context-sensitive language}

Digressing from our main topic, we use Corollary~\ref{new35} to 
establish a language-theoretic property of primitive words.

Recall that a {\em context-sensitive A-grammar} is a triple $\G =
(V,P,S)$ where $V$ is a finite set containing $A$, $S$ is an element
of $V$ that is not in $A$ and $P$ is the {\em set of rules} of the
grammar: a finite set of pairs $(\ell,r) \in V^+ \times V^+$ such that
$$\ell\not\in A^+ \textrm{ and }|\ell| \leq |r|.$$
For all $x,y \in V^+$, we write $x \Rw y$ if there exist $u,v \in V^*$
and $(\ell,r) \in P$ such that $x = u\ell v$ and $y = urv$.  We denote
by $\dmapright{*}$ the transitive and reflexive closure of $\Rw$.  The
language \textit{generated by} $\G$ is
$$L(\G) = \{ w\in A^+ \mid S \; \dmapright{*} \; w \}.$$
A language $L \subseteq A^+$ is said to be {\em context-sensitive} if
it is generated by some context-sensitive $A$-grammar.  As usual, a
language $L\subseteq A^*$ is called {\em context-sensitive} if $L \cap
A^+$ is context-sensitive.

The right and left quotients of a language $L$ by a word $u$ are
defined by
$$u\backslash L = \{ x \in A^* \mid ux \in L\}, \quad L/u = \{ x \in
A^* \mid xu \in L\}.$$  

\bl\label{closure csl}
The class of context-sensitive languages is closed under union,
intersection, concatenation, right and left quotient by a word, $1$-free
substitutions and inverse morphisms.
%
\el

\proof
Closure under union, intersection, concatenation, $1$-free substitutions,
and inverse homomorphisms 
is well-known
\cite[Exercise 9.10]{HU}.  In particular, the family of
context-sensitive languages forms a \textit{trio} \cite[Section
11.1]{HU} and as such, it is closed under \textit{limited erasing}
\cite[Lemma 11.2]{HU}.  By definition, this means that if $k \ge 1$,
$L$ is context-sensitive and $\p$ is a morphism such that $\p(v) \ne
1$ for each $u \in L$ and each factor $v$ of $u$ of length greater
than $k$, then $\p(L)$ is context-sensitive as well.

For the quotients, it suffices to consider letters, hence let $L
\subseteq A^*$, $a\in A$ and $\$\not\in A$.  Let $\sigma$ 
be the substitution that maps $a$ to $\sigma(a) = \{a,\$\}$ and which
fixes every other letter of $A$.  Let also $\p\colon (A\cup\{\$\})^*
\rightarrow A^*$ be the morphism which fixes every letter of $A$ and
erases $\$$.  Then $a\backslash L = \p(\sigma(L) \cap \$A^*)$ and $L/a
= \p(\sigma(L) \cap A^*\$)$.  Since the $\sigma$-images of the letters
are finite, and hence context-sensitive, the languages $\sigma(L) \cap
\$A^*$ and $\sigma(L) \cap A^*\$$ are context-sensitive; moreover $\p$
exhibits limited erasing on these languages, so $a\backslash L$ and
$L/a$ are context-sensitive as well.
\qed

\bp
\label{cslhom}
Let $A$ be a finite alphabet and let $\Gamma$ be a finite set of
endomorphisms of $A^+$. For every $u\in A^+$, $\Gamma^*(u)$ is a
context-sensitive language.
\ep

\proof
Take $b \notin A$. We define a context-sensitive $(A\cup\{ b\})$-grammar $\G =
(V,P,S)$ by $V = A \cup \{ R,S,T\} \cup \{ F_{\p} \mid \p \in \Gamma \}$
and
\begin{align*}
    P = \{ &S \to bF_{\p}uR,\; S \to bub^2,\; F_{\p}a \to \p(a)F_{\p},
    \; F_{\p}R \to TR,\\
    &F_{\p}R \to b^2, \; aT \to Ta,\; bT \to bF_{\p} ; \; a \in A,\; \p
    \in \Gamma\ \}. 
\end{align*}
We show that $L(\G) = b\Gamma^*(u)b^2$.

Clearly, $F_{\p}v \; \dmapright{*}\; \p(v)F_{\p}$ for all $\p \in \Gamma$ and
$v \in A^*$ and so 
$$bvTR \;\dmapright{*}\; bTvR \Rw bF_{\p}vR \;\dmapright{*}\; b\p(v)F_{\p}R
\Rw b\p(v)TR.$$ 
Since $S \Rw bF_{\p}uR \;\dmapright{*}\; b\p(u)F_{\p}R
\Rw b\p(u)TR$ for
every $\p \in \Gamma$, it follows that $S\; \dmapright{*}\;
b\theta(u)F_{\p}R \Rw b\theta(u)b^2$ for
every $\theta \in \Gamma^+$. Together with $S \Rw bub^2$, this yields
$b\Gamma^*(u)b^2 \subseteq L(\G)$.

To prove the opposite inclusion, let
$$Z = \big\{ S \big\} \cup \big\{ bxyb^2, bxTyR, b\p(x)F_{\p}yR \mid
xy \in \Gamma^*(u)\big\}.$$
Then $Z$ is closed under $\Rw$.  That is: if $X \in Z$ and $X \Rw Y$,
then $Y \in Z$.

Since $S \in Z$, it follows that $L(\G) \subseteq Z \cap A^* =
b\Gamma^*(u)b^2$ and so $L(\G) = b\Gamma^*(u)b^2$.  Thus
$b\Gamma^*(u)b^2$ is context-sensitive and by Lemma~\ref{closure csl},
$\Gamma^*(u) = b\backslash (b\Gamma^*(u)b^2) / b^2$ is context-sensitive
as well.
\qed

\bt
\label{csl}
$\oo{P_2}$ is a context-sensitive language.
\et

\proof
Since the class of context-sensitive languages is closed under union
(Lemma~\ref{closure csl}), it follows from Proposition \ref{new34}(ii)
and Proposition \ref{cslhom} that $P_2 \cap \{ a,b\}^+ = \oo{P_2} \cap
\{ a,b\}^+$ is context-sensitive.  Moreover, Proposition
\ref{new34}(i) shows that $P_2 = \Lambda\Psi(P_2 \cap \{ a,b\}^+) =
\Psi\Lambda(P_2 \cap \{ a,b\}^+)$.  Since $\Psi$ is finite, we need
only prove that each $\oo{\psi\Lambda(P_2 \cap \{ a,b\}^+)}$, $\psi\in
\Psi$, is context-sensitive.

Notice that, for each $\psi\in \Psi$ and each word $w$, $\oo{\psi(w)} =
\psi(\oo{w})$.  By Lemma~\ref{closure csl} again, we need only to
prove that $\oo{\Lambda(P_2 \cap \{ a,b\}^+)}$ is context-sensitive.

Let $w\in R_2$ and $p\in P_2 \cap \{ a,b\}^+$.  If $wpw\inv$ is not
reduced, then one of $wp$ and $pw\inv$ is not reduced.  In the first
case, let $q$ be the longest prefix of $p$ such that $q\inv$ is a
suffix of $w$, say $p = qr$ and $w = vq\inv$.  Then $\oo{wpw\inv} =
\oo{vq\inv qrqv\inv} = \oo{vrqv\inv}$.  The second case (if $wp$ is
reduced but $pw\inv$ is not) is treated similarly.  Iterating this
reasoning, we find that $\oo{wpw\inv} = vp'v\inv$, where $v$ is a
prefix of $w$ and $p'$ is a cyclic shift of the word $p$ -- that is,
there are words $q,r$ such that $p = qr$ and $p' = rq$.

Since $P_2 \cap \{ a,b\}^+$ is closed under taking cyclic shifts, it 
follows that $\oo{\Lambda(P_2 \cap \{ a,b\}^+)}$ is the set of 
reduced words of the form $vpv\inv$ with $p\in P_2 \cap \{ a,b\}^+$.

Thus, if $\G = (V,P,S)$ is a context-sensitive $A$-grammar generating
$P_2 \cap \{ a,b\}^+$, then $\oo{\Lambda(P_2 \cap \{ a,b\}^+)} =
L(\G') \cap R_2$, where $\G' = (V',P',S')$ is the context-sensitive
$A$-grammar given by $S'\not\in V$, $V' = \{S'\} \cup V$ and $P' = P
\cup \{ S' \to S\} \cup \{S' \to cS'c\inv; \; c \in A_2 \cup A_2\inv
\}$.  In view of the closure properties in Lemma~\ref{closure csl},
$\oo{\Lambda(P_2 \cap \{ a,b\}^+)}$ is context-sensitive, and hence so
is $P_2$ .
\qed

This result cannot be improved to the next level of Chomsky's
hierarchy:

\bp
\label{nocfl}
$\oo{P_2}$ is not a context-free language.
\ep

\proof
We show that $P_2 \cap ab^+ab^+ab^+$  is not a 
context-free language. Since the class of context-free languages is
closed under intersection with rational languages, it shows that $P_2$
is not context-free either.

It follows easily from Proposition \ref{new34}(i) that $P_2 \cap 
ab^*ab^*ab^*$ is equal to
\beq
\label{nocfl1}
\left\{ ab^mab^nab^{k} \mid  m,n,k \in \N,\; 
\max(m,n,k) = \min(m,n,k)+1\right\}.
\eeq
It is now a classical exercise to show that $P_2 \cap ab^+ab^+ab^+$ is
not context-free since it fails the Pumping Lemma for context-free
languages \cite[Section 6.1]{HU}.
\qed

\section{Singularities, bridges and automorphisms in $\Sigma$}

We now discuss the evolution of the Stallings automaton of a subgroup
$H$ under the iterated action of the automorphisms in $\Sigma$.  It is
well-known that the automata $\A(\p(H))$ may grow unboundedly as the
length of $\p$ (as a product of elements of $\Sigma$) grows.  But in
the context of the mixed orbit problem with respect to the
automorphisms in $\Sigma^*$, we are only interested in the possibility
of reading a $u$-labeled loop (where $u$ is a fixed word) in
$\A(\p(H))$: if the growth of the automata results in long stretches
without branchpoints, then this growth does not affect the membership
of $u$ in $\p(H)$ after a certain point.

Indeed, we show that the fragments of the $\A(\p(H))$
($\p\in\Sigma^*$) that could conceivably allow the reading of a
$u$-loop take only finitely many values -- and these fragments (which
we call truncated automata) can be organised as the states of an
automaton on alphabet $\Sigma$.  The mixed orbit problem with respect
to $\Sigma^*$ then reduces to deciding whether this automaton accepts
a non-empty language.

We now get into the technical considerations that give substance to
this overview of our method. Given $H \lefg F_2$, we say that a
state $q$ of $\A(H)$ is
\bi
\item
a {\em source} if $q\cdot a,q\cdot b \neq \emptyset$, \hspace{2cm}
$\mapleft{a} q 
\mapright{b}$  
\item
a {\em sink} if $q\cdot a\inv,q\cdot b\inv \neq \emptyset$. \hspace{2cm}
$\mapright{a} q 
\mapleft{b}$ 
\ei
Note that a source may have incoming edges and a sink may have
outgoing edges. We use the general term {\em singularities} to refer
to both sources 
and sinks and we denote by $\sing(H)$ the set of all singularities of
$\A(H)$ plus the origin.

If we emphasize the vertices of $\sing(H)$ in $\A(H)$, it is immediate
that $\A(H)$ can be described as the union of {\em positive paths},
i.e. paths with label in $(a\cup b)^+$, between the vertices of
$\sing(H)$, and these positive paths do not intersect each other
except at $\sing(H)$.  We call such paths {\em bridges}.  Note that
every positive path whose internal states are not singularities can be
extended into a uniquely determined bridge.

\subsection{Bridges  in $\A(H)$}

The next two results are easily verified.

\bfa
\label{factuno}
The automaton $\A(\p_{b\inv,a\inv}(H))$ has the same vertex set as
$\A(H)$, edges are reverted and labels changed.  In particular,
sources and sinks are exchanged.  If $\beta$ is a bridge in $\A(H)$,
$\beta = p \mapright{w} q$, then there is a bridge of equal length $q
\mapright{} p$ in $\A(\p_{b\inv,a\inv}(H))$, labeled
$\p_{b\inv,a\inv}(w\inv)$, which we denote by
$\p_{b\inv,a\inv}(\beta)$.
\efa

\bfa
\label{leix}
The automaton $\A(\p_{b,a}(H))$ has the same vertex set as $\A(H)$ and
labels are exchanged.  Sources and sinks remain the same.  If $\beta$
is a bridge in $\A(H)$, $\beta = p \mapright{w} q$, then there is a
bridge of equal length $p \mapright{} q$ in $\A(\p_{b,a}(H))$, labeled
$\p_{b,a}(w)$, which we denote by $\p_{b,a}(\beta)$.
\efa

\noindent Dealing with $\p_{a,ba}$ is naturally a little more complex.
However, as we will see in the next two statements, the foldings
implied in computing $\A(\p_{a,ba}(H))$ are very local: they can be
performed in a single round of independent foldings.  Moreover,
sources in $\A(\p_{a,ba}(H))$ were already sources in $\A(H)$ and
sinks in $\A(\p_{a,ba}(H))$ are at distance 1 of sinks in $\A(H)$.

\bfa
\label{factdos}
The automaton $\A(\p_{a,ba}(H))$ is obtained from
$\A(H)$ by the  
following 3 steps:
\bi
\item[(S1)]
If $p \mapright{b} q$ is an edge of $\A(H)$ and $q$ is not a sink, we
replace that edge by  a path $p \mapright{b} \bullet \mapright{a} q$,
adding a new intermediate vertex for each such edge.
\item[(S2)]
If $p \mapright{b} q \mapleft{a} r$ is a sink in $\A(H)$, we replace
this configuration by
$$\xymatrix{
p \ar@/^/[rr]^{b}& q  &r  \ar@/^/[l]^{a}}$$
\item[(S3)] We iteratively
  remove all the vertices of degree 1 different from the origin.
\ei
\efa

\proof 
Following \cite[Subsection 1.2]{RVW}, the automaton $\A(\p(H))$ may be
obtained from $\A(H)$ in three steps: 
\bi
\item[(1)] We replace each edge labelled by $b$ by a path labelled
  $ba$ (introducing a new intermediate vertex for each such edge),
  producing a dual automaton $\B$.
\item[(2)] We execute the complete folding of $\B$.
\item[(3)] We successively
  remove all the vertices of degree 1 different from the origin.
\ei
How much folding is involved in the process? Let us consider
the first level of folding, i.e. those pairs of edges that can be
immediately identified in $\B$. 
\bi
\item
There are no $b$-edges involved in the first level of folding: indeed,
the $b$-edges keep their origin when we go from $\A(H)$ to $\B$, and
their target is always a new vertex where folding cannot take place.
\item
If we have a sink $p \mapright{b} q \mapleft{a} r$ in $\A(H)$, we get
$$p \mapright{b} \bullet \mapright{a} q \mapleft{a} r$$
in $\B$ and
therefore an instance of  first level folding, yielding
$$\xymatrix{
p \ar@/^/[rr]^{b}& q  &r \ar@/^/[l]^{a} }$$
\item
These are the only instances of  first level folding: we cannot fold
two ``new'' $a$-edges $\mapright{a} q \mapleft{a}$ in $\B$ since that
would imply the existence of two $b$-edges $\mapright{b} q
\mapleft{b}$ in $\A(H)$.
\ei
Let $\C$ denote the automaton obtained by performing all the instances
of first level folding in $\B$. It follows from the above remarks that
$\C$ can be obtained from $\A(H)$ by application of (S1) and (S2).

We actually need no second level of folding because $\C$ is already
deterministic. Indeed, it is clear from (S1) and (S2) that
configurations such as
$\mapleft{a} q \mapright{a}$ or $\mapleft{b} q \mapright{b}$ cannot
occur in $\C$. 

Suppose that $\mapright{b} q \mapleft{b}$ does
occur. Then both edges must have been obtained through (S2) 
and the 
origin of these edges is the vertex $q\cdot(ab\inv)$ in $\A(H)$, a 
contradiction.

Finally, suppose that $\mapright{a} q \mapleft{a}$ does
occur. At least one of these edges must have been obtained through (S1), but
not both, otherwise we would have a configuration $\mapright{b} q
\mapleft{b}$ in $\A(H)$. But then we would have a configuration
$\mapright{a} q \mapleft{b}$ in $\A(H)$ and $q$ would be a sink,
contradicting the application of (S1). Thus $\C$ is deterministic and
so $\A(\p(H))$ is obtained from $\A(H)$ by successive application of
(S1), (S2) and (S3).
\qed

\bfa
\label{facttres}
\bi
\item[(i)]
When applying $\p_{a,ba}$, a state of $\A(H)$ is trimmed in step
(S3) if and only if it is a sink  
of $\A(H)$ without outgoing edges. Moreover, no consecutive states
can be trimmed.
\item[(ii)]
The sources of $\A(\p_{a,ba}(H))$ are precisely the sources $p$ of
$\A(H)$ such that $p\cdot a$ is not a sink or has outgoing edges in
$\A(H)$. 
\item[(iii)]
The sinks of $\A(\p_{a,ba}(H))$ are precisely the states $p$ of
$\A(H)$ with incoming edges such that $p\cdot a$ is a sink of $\A(H)$.
\ei
\efa

\proof
(i) The origin cannot be trimmed and the number of outgoing edges
never decreases, so the only possible candidates to (S3) are the
states that see a decrease in their number of incoming edges, which
are precisely the sinks of $\A(H)$.  Their fate will then depend on
the previous existence of some outgoing edge.  Note that $\A(H)$
cannot possess two consecutive sinks with no outgoing edges, hence the
trimming of a vertex will not be followed by the trimming of any of
its neighbours.

(ii) Since outgoing edges can be at most redirected through (S1) and
(S2), it is clear that every source $p$ of $\A(\p_{a,ba}(H))$ must be
a source of
$\A(H)$. Thus everything will depend on $p\cdot a$ being trimmed or not,
and part (i) yields the claim. 

(iii) No new intermediate vertex obtained through (S1) can become a
sink, and any sink of $\A(H)$ will not remain such after application
of (S2).  Thus the only remaining candidates are the non-sinks of
$\A(H)$ that see an increase of their number of incoming edges, which are
precisely those of the form $q\cdot a\inv$, where $q$ is a sink of
$\A(H)$.  Clearly, to have two distinct incoming edges in
$\A(\p_{a,ba}(H))$, $p = q\cdot a\inv$ must have at least one incoming
edge in $\A(H)$.  In such a case, it is easy to check that after
(S1)/(S2), $p$ has indeed become a sink of $\A(\p_{a,ba}(H))$.  We
remark also that the subsequent trimming by (S3) does not affect the
presence of singularities.
\qed

\bfa
\label{factfour}
Let $\beta = p \mapright{w} q$ be a bridge in $\A(H)$ of length at
least 2, and let $w = w'cd$ where $c,d \in A$.
\bi
\item[(i)]
$\A(\p_{a,ba}(H))$ has a positive path $p \vvlongmapright{\p_{a,ba}(w'c)}
s$, which extends  
to a uniquely determined bridge, denoted by $\p_{a,ba}(\beta)$.
\item[(ii)]
$|\p_{a,ba}(\beta)| \ge |\beta| - 1$, and we have
$|\p_{a,ba}(\beta)| =  
|\beta| - 1$ exactly if $w \in a^+$, $p$ is a source or the
origin in $\A(H)$,  
and $q$ is a sink in $\A(H)$.
\ei
\efa

\proof
Write
$\beta = p 
\mapright{w'} r \mapright{c} s \mapright{d} q$.

(i) By Fact \ref{facttres}, no state of the path $p
\vvlongmapright{\p_{a,ba}(w'c)} s$ risks trimming.  Hence it suffices
to check that no internal state of this path can become a singularity.
This follows easily from Fact \ref{facttres} (ii) and (iii).

(ii) The inequality $|\p_{a,ba}(\beta)| \ge |\beta| - 1$ follows at once
from part (i). It follows also that $|\p_{a,ba}(\beta)| =  
|\beta| - 1$ if and only if $w'c
\in a^+$ (otherwise $|\p_{a,ba}(\beta)| \geq
|\p_{a,ba}(w'c)| \; > \; |w'c| = |\beta| - 1$) and
$p,s \in \sing(\p_{a,ba}(H))$. Thus we assume that $w'c
\in a^+$. 

Clearly, if $p$ is the origin, it must remain so.  If $p$ is a
source, it follows from Fact \ref{facttres} (ii) that $p$ remains a
source (since $p\cdot a$ is not a sink in $\A(H)$).  Finally, if $p$
is a sink, it will no longer be a singularity in $\A(\p_{a,ba}(H))$ by
Fact \ref{facttres} (iii).  Therefore $p \in \sing(\p_{a,ba}(H))$ if
and only if it is a source or the origin in $\A(H)$.

Similarly, $q$ can never become the origin or a source.
Since $q$ has incoming edges in $\A(H)$, it follows from
 Fact \ref{facttres}(iii) that $s$ becomes a sink in
$\A(\p_{a,ba}(H))$ if and only if $s\cdot a$ is a sink in
$\A(H)$. Since the unique outgoing edge of $s$ in $\A(H)$ has label
$d$, then $s \in
\sing(\p_{a,ba}(H))$ if and only if $d = a$ and $q$ is a sink in
$\A(H)$.
\qed

\subsection{Homogeneous cycles and cycle-free paths}

Let $\sigma(H) = \max(1, \source(H)+\sink(H))$, where $\source(H)$
(resp.  $\sink(H)$) is the number of sources (resp.  sinks) of
$\A(H)$.  We call $\sigma(H)$ the {\em number of singularities} of
$\A(H)$.  Note that a vertex may be a source and a sink, and in that
case, it contributes twice to $\sigma(H)$.

We say that a path $p \mapright{w} r$ is {\em homogeneous} if $w \in
R_a \cup R_b$, and it is \textit{special homogeneous} if, in addition,
it starts at a source or the origin, and it ends at a sink or the
origin. Let $\hc(\A)$ (resp. $\hcfp(\A)$, $\shcfp(\A)$) be the maximum 
length of a homogeneous cycle (resp. homogeneous cycle-free path, 
special homogeneous cycle-free path) in automaton $\A$.

Given $H \leq_{f.g.} F_2$, we define
\begin{align*}
    \delta_0(H) &= \max (\sigma(H), \hc(\A(H)),\\
    \delta(H) &= \max (\delta_0(H), \hcfp(\A(H)),\\
    \zeta(H) &= \max (\delta_0(H), \shcfp(\A(H)).
\end{align*}
We record the following inequalities.

\begin{Lemma}\label{b paths}
    Let $H \leq_{f.g.} F_2$.  Every cycle or a cycle-free path labeled
    $b^k$ in $\A(\p_{a,ba}(H))$ satisfies $k \le \sigma(H)$.
\end{Lemma}

\proof
Let us first assume that $\alpha = p\mapright{b^k} q$ is a cycle-free path, say
$$p = q_0 \mapright{b} q_1 \mapright{b} \ldots \mapright{b} q_k = q.$$
Since any
$b$-edge obtained through (S1) must be followed only by an $a$-edge
(see Fact~\ref{factdos}), only the last edge $q_{k-1} \mapright{b}
q_k$ may be obtained through (S1), and the other edges arise from
applications of (S2).  Thus there exist edges in $\A(H)$ (represented
through discontinuous lines) of the form
$$\xymatrix{ &p_1 &p_2 && p_{k-2} & p_{k-1} &\\
q_0 \ar@{-->}[ur]^b \ar[r]_b & q_1 \ar@{-->}[ur]^b \ar@{-->}[u]_a \ar[r]_b & q_2
\ar@{-->}[u]_a & \ldots & q_{n-2} \ar[r]_b \ar@{-->}[u]_a
\ar@{-->}[ur]^b & q_{k-1} 
\ar@{-->}[u]_a  
\ar[r]_b & q_{k} }$$
In particular, the vertices $p_1,\ldots, p_{k-1}$ are distinct sinks
in $\A(H)$, and the vertices $q_1,\ldots, q_{k-1}$ are distinct
sources in $\A(H)$.  Therefore $2k-2 \le \sigma(H)$ and hence $k \le
\sigma(H)$.

If $\alpha$ is a cycle, then not even the last edge of $\alpha$ arises
from an application of (S1), and the same reasoning shows that $2k \le
\sigma(H)$, so $k\le \sigma(H)$.
\qed

\bl
\label{deltazero}
Let $H \leq_{f.g.} F_2$ and $\p \in \Sigma$.  Then
\begin{align*}
    \sigma(\p(H)) &\leq \sigma(H),\\
    \delta_0(\p(H)) &\leq \delta_0(H),\\
    \zeta(\p(H)) &\leq \zeta(H).
\end{align*}
\el

\proof
The first inequality is a direct consequence of Facts \ref{factuno},
\ref{leix} and \ref{facttres}.

By Facts \ref{factuno} and \ref{leix}, the other inequalities are
trivial if $\p = \p_{b,a}$ or $\p_{b\inv,a\inv}$.  We now assume that
$\p = \p_{a,ba}$.  Since $\sigma(\p(H)) \leq \sigma(H)$, we only need
to show that the maximum length of a homogeneous cycle (resp.
cycle-free special homogeneous) path $\alpha = p \mapright{} q$ in
$\A(\p(H))$ ($p = q$ in the case of a cycle) is at most equal to
$\delta_0(H)$ (resp.  $\zeta(H)$).

If the label of $\alpha$ is $b^k$, then Lemma~\ref{b paths} shows 
that $k\le \sigma(H)$, so $k \le \delta_0(H) \le \zeta(H)$.

Suppose now that the label of $\alpha$ is $a^k$.  In view of Fact
\ref{factdos}, none of its edges was obtained trough (S1): indeed the
$a$-edge in $\bullet \mapright{b} \bullet \mapright{a} \bullet$
produced by (S1) cannot occur in a homogeneous cycle, nor in a
homogeneous path unless it is its first edge.  But its initial vertex
is not a singularity, so this edge cannot occur in a special
homogeneous path.  Hence the path $\alpha$ already existed in 
$\A(H)$. If $\alpha$ is a cycle, then $k \le \delta_0(H)$.

If instead $\alpha$ is a special homogeneous cycle-free path, then
Fact \ref{facttres} (ii) shows that $p$ is either the origin or a
source in $\A(H)$.  If $q$ is the origin, we immediately get $k \leq
\zeta(H)$.  If instead $q$ is a sink in $\A(\p(H))$, then $s = q\cdot
a$ is a sink of $\A(H)$ by Fact \ref{facttres} (iii), and we have a
path
$$\alpha' \enspace = \enspace p \vlongmapright{a^k} q \longmapright{a}
s$$
in $\A(H)$.  If $\alpha'$ is cycle-free, then $k < k+1 \leq \zeta(H)$.
If, on the contrary, $\alpha'$ is not cycle-free, then $s$ is the only
repetition since the length $k$ prefix of $\alpha'$, namely $\alpha$,
is cycle-free.  If $s \neq p$, then $q$ would also be a repetition
since $\A(H)$ is an inverse automaton.  Therefore $s = p$, so
$\alpha'$ is a homogeneous cycle in $\A(H)$ and hence $k < k+1 \leq
\delta_0(H) \leq \zeta(H)$.  This concludes the proof.
\qed

\begin{Remark}\rm
    Note that it is not the case that $\delta(\p(H)) \leq \delta(H)$
    always holds when $\p\in\Sigma$: see the case where $H = \langle
    ba \rangle$ and $\p = \p_{a,ba}$.
\end{Remark}

\subsection{Truncated automata}\label{sec: truncated}

Given $H \leq_{f.g.} F_2$, we consider the {\em geodesic metric} $d$
defined on the vertex set of $\A(H)$ by taking $d(u,v)$ to be the
length of the shortest path connecting $u$ and $v$.  Since $\A(H)$ is
inverse, it is irrelevant to consider directed or undirected paths.
As usual, we have
$$d(u,\sing(H)) = \min\{ d(u,v) \mid v \in \sing(H) \}.$$
Given $t > 0$, the $t$-{\em truncation} of $\A(H)$, denoted by
$\A_t(H)$, is the automaton obtained by removing from $\A(H)$ all
vertices $u$ such that $d(u,\sing(H)) > t$ and their adjacent
edges. Note that this automaton does not need to be connected.

We first observe that if $\beta$ is a bridge which is long enough to
be affected by the $t$-truncation of $\A(H)$, then for each
$\p\in\Sigma$, $\p(\beta)$ is affected by the $t$-truncation of
$\A(\p(H))$ as well.

\bp
\label{mainlemma}
Let $\p \in \Sigma$, $H \leq_{f.g.} F_2$ and $K \in \Sigma^*(H)$. 
If $\beta$ is a bridge in $\A(K)$ and $|\beta|  > \zeta(H)$,
  then $|\p(\beta)| \geq |\beta|$.
\ep

\proof
The result is trivial if $\p = \p_{b\inv,a\inv}$ or $\p = \p_{b,a}$
since in those cases, $|\p(\beta)| = |\beta|$ (Facts \ref{factuno} and
\ref{leix}).  We now assume that $\p = \p_{a,ba}$.

By Fact \ref{factfour}, if $|\p(\beta)| < |\beta|$, then $\beta = p
\mapright{a^k} q$, where $k > \zeta(H)$, $p$ is a source or the
origin in $\A(K)$, and $q$ is a sink of $\A(K)$.  In particular, 
$\beta$ is a special homogeneous cycle-free path, so that $|\beta| 
\le \zeta(K)$.

Since $K \in \Sigma^*H$, Lemma~\ref{deltazero} shows that $\zeta(K)
\le \zeta(H)$, a contradiction.
\qed

\bt
\label{newmainlemma}
Let $\p \in \Sigma$, $H \leq_{f.g.} F_2$, $t > \frac{1}{2}\zeta(H)$
and $K,K' \in \Sigma^*(H)$.  Then
$$\A_t(K) = \A_t(K') \enspace\Longrightarrow\enspace
\A_t(\p(K)) = \A_t(\p(K')).$$
\et

\proof
As in several previous proofs, the result is trivial if $\p =
\p_{b,a}$ or $\p_{b\inv,a\inv}$, and we may assume that $\p =
\p_{a,ba}$.

By Proposition~\ref{mainlemma}, we know that, once the length of a
bridge reaches the threshold $\zeta(H)+1$, it can only get longer.
Since $t > \frac{1}{2}\zeta(H)$, $t$-truncation affects only bridges
of length at least 
$\zeta(H)+1$.  We must therefore discuss the truncation mechanism for
such long bridges.

Assume that $\beta = p \mapright{w} q$ is a bridge in $\A(\mu(H))$
$(\mu \in \Sigma^*)$ with $|w| \geq 2t+1$.  Then we may write $w =
uzv$ with $|u| = |v| = t$.  By Proposition~\ref{mainlemma}, the label
of $\p(\beta)$ is of the form $u'z'v'$ with $|u'| = |v'| = t$ and
$|z'| \geq |z|$.  We only need to prove that $u'$ and $v'$ depend only
on $\A_t(\mu(H))$ and are therefore independent from $z$.

In view of Fact \ref{facttres}, it is clear that $u'$ depends only on
$\A_t(\mu(H))$ (remember that $w = uzv$ is a positive word and
singularities cannot {\em move forward} along a positive path). The
nontrivial case is of course the case of $q$ being a sink in
$\A(\mu(H))$, since by Fact \ref{facttres} (iii) a sink can actually be
transferred to the preceding state along a positive path. We claim
that even in this case $v'$ is independent from $z$. 

Indeed, assume first that $b$ occurs in $v$.  Then $|\p(v)| > |v|$
provides enough compensation for the sink moving backwards one
position.  Hence we may assume that $v = a^{t}$.  We claim that $v' =
a^t$ as well, independently from $z$.  Suppose not.  Since we are
assuming that the sink has moved from $q$ to its predecessor, and
$\p(a^{t-1}) = a^{t-1}$, it follows that $v' = ba^{t-1}$.  Hence $b$
occurs in $w$.  Write $w = xba^m$.  Since $\p(ba^m) = ba^{m+1}$, and
taking into account the mobile sink, we obtain by comparison $ba^{m} =
ba^{t-1}$ and so $m = t-1$, a contradiction, since $a^t$ is a suffix
of $w$.  Therefore $v' = a^t$ and so is independent from $z$ as
required.  \qed

\bc
\label{newmainlemmaA}
Let $H \leq_{f.g.} F_2$ and $t > \frac{1}{2} \zeta(H)$. Then the
set
$$\X(t,H) = \{ \A_t(K) \mid K \in \Sigma^*(H) \}$$
is finite and effectively constructible.
\ec

\proof
By Lemma \ref{deltazero}, every automaton $\A(K)$, $K\in\Sigma^*(H)$, has 
at most $\sigma(H)$ singularities. By definition of a $t$-truncation, 
every state in $\A_t(K)$ is at distance at most $t$ from a 
singularity, and hence the size of $\A_t(K)$ is bounded. Thus $\X(t,H)$ is finite.

The proof of Theorem \ref{newmainlemma} provides a straighforward
algorithm to compute all its elements.  Indeed, all we need is to
compute the finite sets
$$\X_n(t,H) = \{ \A_t(K) \mid K \in \Sigma^n(H) \}$$
until reaching 
\beq
\label{newmainlemma1}
\X_{n+1}(t,H) \subseteq \bigcup_{i=0}^n \X_i(t,H),
\eeq
which must occur eventually since $\X(t,H) = \cup_{i\geq 0} \X_i(t,H)$
is finite. Why does (\ref{newmainlemma1}) imply $\X(t,H) = \cup_{i\geq
  0}^n \X_i(t,H)$? Suppose that $\B \in \X_m(t,H) \setminus (\cup_{i\geq
  0}^n \X_i(t,H))$ with $m$ minimal, say
$\B = \A_t(\p(K))$ with $K \in \Sigma^{m-1}(H)$ and $\p \in \Sigma$. By
minimality of $m$, we have $\A_t(K) \in \cup_{i\geq
  0}^n \X_i(t,H)$. Thus $\A_t(K) = \A_t(K')$ for some $K' \in
\cup_{i=0}^n \Sigma^i(H)$. Now Theorem \ref{newmainlemma} yields
$$\B = \A_t(\p(K)) = \A_t(\p(K')) \in \bigcup_{i=0}^{n+1} \X_i(t,H) =
\bigcup_{i=0}^n \X_i(t,H),$$
a contradiction. Therefore $\X(t,H) = \cup_{i\geq
  0}^n \X_i(H)$ as claimed.
\qed

\section{Back to orbit problems in $F_2$}
\label{sec: orbit pbms}

We saw in Section~\ref{sec: orbit pbms equation} that it is decidable
whether a given element $u\in F_2$ has an automorphic image in a given
rational subset of $F_2$, and in particular in a given finitely
generated subgroup of $F_2$ (Corollary~\ref{dorbit} above).  We use
truncated automata to give a different proof of this result in the
finitely generated subgroup case.  We also prove the decidability of
mixed orbit problems under the action of certain rational subsets of
$\Sigma^*$.

\subsection{Some mixed orbit problems}
\label{sec: rational autom}

The archetypal result in this section is the solution of the following
mixed orbit problem.

\begin{Proposition}\label{prop: pb}
    Let $u\in F_2$ an $H\lefg F_2$.  The set of automorphisms $\phi\in
    \Sigma^*$ such that $u \in \p(H)$ (resp.  a conjugate of $u$ lies
    in $\p(H)$), is an effectively constructible rational subset of
    $\Sigma^*$.
\end{Proposition}

\proof
Let $t > \max(\frac12\zeta(H),\frac12|u|)$: if $\p\in \Sigma^*$, then
$u\in \p(H)$ if and only if $u$ labels a loop at the origin in
$\A(\p(H))$.  Note that this is the case if and only if $u$ labels a
loop at the origin in $\A_t(\p(H))$.

We now view $\Sigma$ as a finite alphabet (besides being a subset of
$\aut F_2$) and we consider the $\Sigma$-transition system $\B_t(H)$
defined as follows.  (A \textit{$\Sigma$-transition system} is defined
like a $\Sigma$-automaton, omitting the specification of the initial
and terminal states.)  The state set of $\B_t(H)$ is $\X(t,H)$ (see
Corollary~\ref{newmainlemmaA}) and its transitions are the triples
$\A_t(K) \mapright{\p} \A_t(\p(K))$, for each $\A_t(K) \in \X(t,H)$
and $\p\in\Sigma$.  Note that $\X(t,H)$ is finite and effectively
constructible by Corollary~\ref{newmainlemmaA} and the transitions of
$\B_t(H)$ are well-defined by Theorem~\ref{newmainlemma}.  Moreover,
this transition system is complete and deterministic by construction
(and so defines recognizable subsets of $\Sigma^*$ as a submonoid of
$\aut F_2$).
It is immediate that if the word $(\p_1,\ldots,\p_n) \in \Sigma^*$
labels a path from $\A_t(K)$ to $\A_t(K')$ in $\B_t(H)$ ($K,K'\in
\Sigma^*(H)$), then $\A_t(K') = \A_t(\p_n\ldots \p_1(K))$.

Now consider the automaton formed by the transition system $\B_t(H)$
with initial state $\A_t(H)$ and terminal states the elements
$\A_t(K)$ of $\X(t,H)$ such that $u$ labels a loop at the origin in
$\A_t(K)$ (i.e. $u\in K$).  The above discussion shows that the
language accepted by this automaton is the set of words
$(\p_1,\ldots,\p_n) \in \Sigma^*$ such that $u \in\p_n\cdots\p_1(H)$. 
Thus the set of all $\p\in \Sigma^*$ such that $u\in \p(H)$ is 
rational and effectively constructible.

Observe that a conjugate of $u$ lies in $\p(H)$ ($\p\in\aut F_2$), if
and only if $u$ labels a loop at the origin in $\A(\lambda_w\p(H))$
for some $w\in F_2$, if and only if $\cc(u)$ labels a loop somewhere
in $\A(\p(H))$.  We now consider the $\Sigma$-transition system
$\B_t(H)$ as above, with the same initial state, and we take as
terminal states the elements $\A_t(K)$ of $\X(t,H)$ such that
$\cc(u)$ labels a loop anywhere in $\A_t(K)$.  The language in
$\Sigma^*$ accepted by the resulting automaton is the set of
$\p\in\Sigma^*$ such that $\p(H)$ contains a conjugate of $u$.
\qed

The same idea --- and the same transition system--- can be used to
algorithmically solve a number of other orbit problems.

\begin{T}\label{many decidability results}
    Let $H, K \lefg F_2$, $u, u_1,\ldots,u_k \in F_2$ and $R \in
    \rat \Sigma^*$.  Then the following problems are decidable:
    \begin{itemize}
	\item[$(1)$] whether $u \in \mu(H)$ for some $\mu \in R$;
	
	\item[$(1')$] whether a conjugate of $u$ lies in $\mu(H)$ for 
	some $\mu \in R$; that is, whether $u \in \mu(H)$ for some 
	$\mu \in \Lambda R$;

	\item[$(2)$] whether $K \subseteq \mu(H)$ for some $\mu \in R$;
	
	\item[$(2')$] whether a conjugate of $K$ is contained in
	$\mu(H)$ for some $\mu \in R$; that is, whether $K \subseteq
	\mu(H)$ for some $\mu \in \Lambda R$;

	\item[$(3)$] whether $K = \mu(H)$ for some $\mu \in R$;
	
	\item[$(3')$] whether a conjugate of $K$ is equal to $\mu(H)$
	for some $\mu \in R$; that is, whether $K = \mu(H)$ for some
	$\mu \in \Lambda R$;

	\item[$(4)$] whether there exist $w_1,\ldots,w_k\in F_2$ such
	that $\lambda_{w_1}(u_1),\ldots,\lambda_{w_k}(u_k) \in \mu(H)$
	for some $\mu \in R$.
    \end{itemize}
    In addition, for each of these problems, the set of morphisms 
    $\mu\in\Sigma^*$ that it defines is rational and effectively 
    constructible. 
\end{T}

\proof
The solutions of Problems $(1)$ and $(1')$ follow from
Proposition~\ref{prop: pb}: the set $X$ of automorphisms $\p\in\Sigma^*$
such that $\p(H)$ contains $u$ (resp.  a conjugate of $u$) is rational
and we can compute a $\Sigma$-automaton recognizing that set.  Since
$\B_t(H)$ is deterministic and complete, we only
have to decide whether $X$ has a non-empty intersection with the
given rational set $R$, a classical decidable result from automata
theory.

The other proofs follow the same pattern, and correspond to variants
of Proposition~\ref{prop: pb}.  Let us consider Problem $(2)$ and let
$u_1,\ldots,u_k$ be generators of $K$.  Then we need to consider the
$\Sigma$-transition system $\B_t(H)$ with $t >
\max(\frac12\zeta(H),\frac12|u_1|,\ldots,\frac12|u_k|)$, and to choose
as terminal states the elements $\A\in \X(t,H)$ such that $\bar
u_1,\ldots,\bar u_k$ label loops at the origin in $\A$.

For Problem $(2')$, we consider a \textit{cyclically reduced}
conjugate $K'$ of $K$, that is, one such that the origin in $\A(K')$
has degree at least 2 (if the origin in $\A(K)$ has degree 1, choose
any vertex $v$ with degree at least 2 as the new origin and let $K'$
be the corresponding conjugate).  Let $u_1,\ldots,u_k$ be generators
of $K'$.  Then a conjugate of $K$ lies in $\mu(H)$ if and only if the
$\bar u_i$ label loops around the same vertex of $\A(\mu(H))$.  Thus
it suffices to choose $t >
\max(\frac12\zeta(H),\frac12|u_1|,\ldots,\frac12|u_k|)$, and to take
as terminal states the elements $\A\in \X(t,H)$ such that $\bar
u_1,\ldots,\bar u_k$ label loops around the same vertex of $\A$.

For Problem $(3)$, we choose again $t >
\max(\frac12\zeta(H),\frac12|u_1|,\ldots,\frac12|u_k|)$, where
$u_1,\ldots,u_k$ are generators of $K$.  In particular, $t$ is large
enough to have $\A_t(K) = \A(K)$, and we choose a single terminal
state, $\A_t(K)$ (if $\A_t(K) \in \X(t,H)$; if that is not the case,
then Problem $(3)$ is decidable, in the negative).  Then we have an
automaton which recognizes the set $L(t,K)$ of all $\mu\in\Sigma^*$
such that $\A_t(\mu(H)) = \A_t(K) = \A(K)$.  Observe now that
truncation creates (pairs of) degree 1 vertices: the automorphisms
$\mu\in L(t,K)$ are such that $\A_t(\mu(H))$ has at most one degree 1
vertex (the origin), and hence $\A_t(\mu(H)) = \A(\mu(H))$.  Thus our
automaton recognizes the set of all $\mu\in\Sigma^*$ such that
$\A(\mu(H)) = \A(K)$, that is, such that $\mu(H) = K$.

For Problem $(3')$, we consider a cyclically reduced conjugate $K'$ of
$K$ and an integer $t$ as in Problem $(2')$.  Again, we have $\A_t(K')
= \A(K')$.  We choose as terminal states the elements of $\X(t,H)$ of
the form $\A_t(\lambda_w(K))$ ($w\in F_2$).  These automata are of one
of the following types:

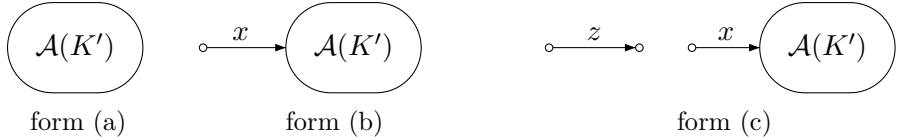
\begin{picture}(53,24)(0,-24)
    \node[Nw=18,Nh=12,Nmr=8](a1)(5,-10){$\A(K')$}
    \node[Nw=18,Nh=12,Nmr=8](a2)(42,-10){$\A(K')$}
    \node[Nw=18,Nh=12,Nmr=8](a3)(105,-10){$\A(K')$}
    \node[Nw=1.0,Nh=1.0,Nmr=0.5](n2)(22.0,-10){}
    \node[Nw=1.0,Nh=1.0,Nmr=0.5](n3)(68.0,-10){}
    \node[Nw=1.0,Nh=1.0,Nmr=0.5](n4)(80.0,-10){}
    \node[Nw=1.0,Nh=1.0,Nmr=0.5](n5)(87.0,-10){}
    \put(-1,-21){\textrm{\small form (a)}}
    \put(33,-21){\textrm{\small form (b)}}
    \put(85,-21){\textrm{\small form (c)}}
    \drawedge[ELpos=25](n2,a2){$x$}
    \drawedge(n3,n4){$z$}
    \drawedge[ELpos=25](n5,a3){$x$}
\end{picture}

\noindent with $|x| \le t$ and $|z| = t$.  As in the discussion of
Problem $(3)$, the existence of a $\mu$-labeled path in $\B_t(H)$ from
$\A_t(H)$ to an automaton of type (a) or (b) shows that $\mu(H)$ is a
conjugate of $K'$, and hence of $K$.  If the path in $\B_t(H)$ ends in
an automaton of type (c), then $\mu(H)$ is a conjugate of $K'$ of the
form $zyxwK'(zyxw)\inv$ or $z\inv yxwK'(z\inv yxw)\inv$ for some $y,w$
such that $zyxw$ or $z\inv yxw$ is reduced.  We then conclude the
proof of the decidability of Problem $(3')$ as usual.

Finally, for Problem $(4)$, we choose $t >
\max(\frac12\zeta(H),\frac12|\cc(u_1)|,\ldots,$ $\frac12|\cc(u_k)|)$ and
we choose as terminal states the elements $\A\in\X(t,H)$ such that
each $\cc(u_i)$ ($i = 1,\ldots,k$) labels a loop at some vertex in $\A$.
\qed

We can also consider finitely many subgroups $H_i$ in (4) and many
other variations.

A simple rewriting of Theorem~\ref{many decidability results} in 
terms of orbit problems (see the introduction) yields the following 
corollary.

\bc
\label{grimc}
Let $H \lefg F_2$, $u \in F_2$ and $R \in \rat \Sigma^*$.  Then
it is decidable whether the orbit of $u$ under the action of $R\inv$
(resp.  $\Lambda R\inv$) meets $H$.

If in addition $K \lefg F_2$, then it is decidable whether $H$
contains an element of the orbit of $K$ under the
action of $R\inv$ or $\Lambda R\inv$; whether $H$ is contained in an 
element of the orbit of $K$ under the action of $R$ or $\Lambda R$; 
and whether $K$ is an element of the orbit of $H$ under the action of 
$R$, $R\inv$, $\Lambda R$ or $\Lambda R\inv$.
\ec

Applying Corollary~\ref{grimc} to the case where $u$ is a letter in 
$A$, we get a statement about primitive elements.

\bc
\label{grimb}
Let $H \lefg F_2$ and $R \in \rat \Sigma^*$.  Then it is
decidable whether $H$ contains a primitive element of the form
$\mu(a)$, $\mu\inv \in R$ (resp. $\mu\inv \in\Lambda R$).
\ec

\begin{Remark}\rm
    Let $S$ be a subset of $R_2$ such that, for each rational set
    $S'$, one can decide whether $S \cap S'$ is empty or not.  Then
    Problems $(1')$, $(2')$ and $(3')$ in Theorem~\ref{many
    decidability results} are decidable even if we restrict the
    conjugating factors to be in $S$, that is, if we replace $\Lambda$
    by $\{\lambda_s \mid s\in S\}$ in the statement of these problems.
    The same restriction can be imposed to $\Lambda$ in the statements
    of Corollaries~\ref{grimc} and~\ref{grimb}.

    Similarly, Problem $(4)$ in Theorem~\ref{many decidability
    results} remains decidable even if we require the $w_i$ to be in
    fixed subsets $S_i$ ($i = 1,\ldots,k$) such that, for each
    rational set $S'$, one can decide whether $S_i \cap S'$ is empty
    or not.
    \enspace\qed
\end{Remark}

\subsection{Orbits under invertible substitutions}\label{sec: invertible}

Invertible substitutions are an interesting special case of the
rational subsets of $\aut F_2$ discussed in Section~\ref{sec: rational
autom}.  This leads to the following statement.

\begin{Corol}
    The problems discussed in Theorem~\ref{many decidability results}
    and Corollaries~\ref{grimc} and~\ref{grimb} are decidable also if
    $R$ is assumed to be a rational subset of $\IS(F_2)$ or
    $\IS(F_2)\inv$.
\end{Corol}

\proof
If $R \in \rat \IS(F_2)$, then $R \in \rat\Sigma^*$ by
Lemma~\ref{lemma: IS and ISinv}, and we simply apply Theorem~\ref{many
decidability results} and Corollaries~\ref{grimc} and~\ref{grimb}.

If $R \in \IS\inv(F_2)$, then $R = \p_{a,b\inv}\ R'\ \p_{a,b\inv}$ for
some $R'\in \rat\Sigma^*$ by Lemma~\ref{lemma: IS and ISinv} ($R'$ is
the set of inverses of the elements of $R$).  Problem $(1)$ in
Theorem~\ref{many decidability results} on instance $u$, $H$ and $R$,
for example, is equivalent to the same problem on instances
$\p_{a,b\inv}(u)$, $\p_{a,b\inv}(H)$ and $R'$, which we know to be
decidable.  The other problems are handled in the same fashion.
\qed

\subsection{Another solution of the mixed orbit problem for $\aut F_2$}\label{sec: proof}

Our proof relies on truncated automata and Theorem \ref{deco}. The key
is to bound the powers of $\p_{a,ba}$ that we need to consider, and is
achieved in view of our previous bound for the length of homogeneous cycles. 

Let $u \in F_2$ and $H \leq_{f.g.} F_2$. We want to show that it is 
decidable whether $\mu(u) \in H$ for some $\mu\in \aut F_2$.
By Theorem \ref{deco}, and since $\Psi\inv = \Psi$, it suffices to
decide whether there exist $w \in F_2$ and $n \geq 0$ such that
one of the
following conditions hold:
\bi
\item
$\lambda_w\p_{a,ba}^n\p_{a\inv,b}(u) \in \Sigma_0^*\Psi(H)$;
\item
$\lambda_w\p_{a,ba}^n\p_{a\inv,b\inv}(u) \in \Sigma_0^*\Psi(H)$.
\ei
Since $\Psi$ is finite, it suffices to be able to decide whether
\beq
\label{dorbit1}
\mbox{there exist $w \in F_2$, $n \ge 0$ and $\mu \in \Sigma_0^*$ such that
  $\lambda_w\p_{a,ba}^n(u) \in \mu(H)$.}
\eeq

We start by considering the case $n = 0$.  By Proposition
\ref{proptec} (i), we may replace $\lambda_w\p_{a,ba}^n$ by
$\p_{a,ba}^n\lambda_w$, so we may assume that $u$ is cyclically
reduced.  And by Proposition \ref{RSS} (iv), our problem further
reduces to asking if one can decide whether
\beq
\label{dorbitD}
\mbox{$u$ labels a loop in $\A(\mu(H))$ for some $\mu\in \Sigma_0^*$.}
\eeq
We note that every loop contains either the origin or a singularity:
if it does not contain the origin, then there is a path from the
origin to a state in the loop, and the first contact between that path
and the loop is a source or a sink.  Now let us fix $t > \max
(\frac12\zeta(H), \frac12|u|)$: then $u$ labels a loop in $\A(\mu(H))$ if and
only if $u$ labels a loop in $\A_t(\mu(H))$.  By the appropriate
variant of Corollary \ref{newmainlemmaA} (where $\Sigma$ is replaced
with $\Sigma_0^*$) we can effectively compute the finite set
$$\X_0(H) = \{ \A_t(K) \mid K \in \Sigma_0^*(H) \}.$$
Thus (\ref{dorbitD}) is decidable, and hence (\ref{dorbit1}) is
decidable for $n=0$.  It is also decidable for any fixed $n$ (applying
the case $n=0$ to $\p_{a,ba}^n(u)$ instead of $u$).

We now consider (\ref{dorbit1}) in its full generality.  If $u \in
R_a$, then we are reduced to the case $n = 0$ since $\p_{a,ba}(u) =
u$.  So we assume that $b$ or $b\inv$ occurs in $u$, and by
conjugation again, we may assume that $u$ starts with $b$ or ends with
$b\inv$ (and not both since $u$ is cyclically reduced).

Let $M$ be the least common multiple of $1, 2, \ldots, \delta_0(H)$.
In order to prove (\ref{dorbit1}), it suffices to show that
\begin{quote}
    \textsl{if there exist $w \in F_2$, $n \geq 0$ and $\mu\in
    \Sigma_0^*$ such that $\lambda_w\p_{a,ba}^n(u) \in \mu(H)$, then
    there exists such a triple $(w,n,\mu)$ with $n < |u| + \max
    (M,\delta(H))$.}
\end{quote}
Since we have proved (\ref{dorbit1}) for bounded $n$, the latter
property is decidable, and hence (\ref{dorbit1}) is decidable in
general.

So we are left with the task of proving this reduced claim.  Let
$(w,n,\mu)$ be such that $\lambda_w\p_{a,ba}^n(u) \in \mu(H)$, with
$n$ minimal, and let us suppose that $n \geq |u| + \max (M,\delta(H))$.

Write $u = a^{i_0}b^{\epsilon_1}a^{i_1} \ldots b^{\epsilon_k}a^{i_k}$
with $k \ge 1$ and $\epsilon_{\ell} = \pm 1$ for every $\ell$.  If
$m\ge 0$, then
$$\p_{a,ba}^m(u) = \p_{a,ba^m}(u) = a^{j_0}b^{\epsilon_1}a^{j_1}
\ldots b^{\epsilon_k}a^{j_k}$$
with
$$j_\ell =
\begin{cases}
    i_\ell + m & \mbox{if $\epsilon_\ell = \epsilon_{\ell+1} = 1$, or
    $\ell = k$ and $\epsilon_k = 1$}\cr
    i_\ell - m & \mbox{if $\epsilon_\ell = \epsilon_{\ell+1} = -1$, or
    $\ell = 0$ and $\epsilon_1 = -1$}\cr
    i_\ell & \mbox{in all other cases.}
\end{cases}$$
Recall that $u$ is cyclically reduced, and that it starts with $b$
($i_0 = 0$ and $\epsilon_1 = 1$) or ends with $b\inv$ ($i_k = 0$ and
$\epsilon_k = -1$).  It follows that $\p_{a,ba^m}(u)$ is cyclically
reduced and that it too starts with $b$ or ends with $b\inv$.

By Proposition \ref{RSS} (iv), $\p_{a,ba}^n(u)$ labels a loop $\alpha$
in $\A(\mu(H))$.  Moreover, we have
$$\p_{a,ba}^n(u) = a^{r_0}b^{\epsilon_1}a^{r_1} \ldots
b^{\epsilon_k}a^{r_k},\quad \p_{a,ba}^{n-M}(u) =
a^{s_0}b^{\epsilon_1}a^{s_1} \ldots b^{\epsilon_k}a^{s_k},$$
with
$$\begin{cases}
    r_\ell = i_\ell + n,\ s_\ell = r_\ell-M  & \mbox{if $\epsilon_\ell = \epsilon_{\ell+1} = 1$, or
    $\ell = k$ and $\epsilon_k = 1$}\cr
    r_\ell = i_\ell - n,\ s_\ell = r_\ell + M & \mbox{if
    $\epsilon_\ell = \epsilon_{\ell+1} = -1$, or $\ell = 0$ and
    $\epsilon_1 = -1$}\cr
    s_\ell = r_\ell = i_\ell & \mbox{in all other cases.}
\end{cases}$$
In the first and second cases, $|r_\ell| > n-|u|\ge
\max(M,\delta(H))$; and in the last case, $|r_\ell| < |u|$.  Thus, for
the indices $\ell$ such that $r_\ell\ne s_\ell$, we have $r_\ell >
\delta(H)$.  We now show that the fragments of the loop $\alpha$
labeled by the factors $a^{r_\ell}$ such that $r_\ell\ne s_\ell$, fail
to be cycle-free in $\A(\mu(H))$.

Recall that $\mu\in \Sigma_0^*$.  If $\mu = \id$ or
$\p_{b\inv,a\inv}$, the result is immediate since $r_\ell > \delta(H)
= \delta(\mu(H))$.  If $\mu = \p_{a,ba}\nu$ with $\nu\in \Sigma_0^*$,
then we can use Proposition~\ref{proptec} (i) to reduce $n$, a
contradiction.  Hence we may assume that $\mu = \p_{b\inv,a\inv}\nu$
with $\nu\in \Sigma_0^*$, $\nu\ne\id$.  Since $\p_{b\inv,a\inv}^2 =
\id$, we may further assume that $\mu = \p_{b\inv,a\inv}\p_{a,ba}\nu'$
with $\nu' \in \Sigma_0^*$.  Then the vertices involved in the
$a^{r_\ell}$-labeled fragment of $\alpha$ form a path in
$\A(\p_{a,ba}\nu'(H))$ labeled $b^{r_\ell}$.  Since $r_\ell >
\delta(H)$, we also have $r_\ell > \sigma(H) \ge \sigma(\nu'(H))$
(Lemma~\ref{deltazero}), and hence this path is not cycle-free by
Lemma~\ref{b paths}.

So, for each $\ell$ such that $r_\ell \ne s_\ell$, the fragment of
$\alpha$ labeled by the factor $a^{r_\ell}$ of $\p_{a,ba}^n(u)$ fails
to be cycle-free, and must be read along a cycle of $\A(\mu(H))$ (in
an inverse automaton, if a homogeneous path contains a cycle, then it
reads entirely along that cycle).

By definition, $M$ is a multiple of the length $c_\ell$ of that cycle.
Now compare $\p_{a,ba}^{n-M}(u)$ and $\p_{a,ba}^n(u)$: wherever the
$a$-factors $a^{r_\ell}$ and $a^{s_\ell}$ are different, their
difference is either $a^M$ or $a^{-M}$, and hence it consists of a
whole number of passages around the length $c_\ell$ cycle.  Therefore
$\p_{a,ba}^{n-M}(u)$ labels a path in $\A(\mu(H))$ as well.  This
contradicts the minimality of $n$ and completes the proof.

\begin{Remark}\rm
    The \textit{a priori} complexity of the algorithms discussed in 
    Section~\ref{sec: orbit pbms} is very high: if $u$ has length at most 
    $n$ and $\A(H)$ has at most $n$ states, then $\sigma(H), 
    \zeta(H) \le n$ and the truncated automata can have exponentially 
    many states. There can  therefore be super-exponentially many 
    truncated automata, forming the states of the transition system 
    $\B_t(H)$ -- in which we must solve a reachability problem 
    (polynomial in the number of states of the transition system).
\end{Remark}

\end{document}